\documentclass[12pt,reqno]{amsart}
\usepackage{epsf, amsmath, amsfonts, amsthm, amssymb}
\newtheorem{lemma}{Lemma}[section]
\newtheorem{theorem}[lemma]{Theorem}
\newtheorem{proposition}[lemma]{Proposition}
\newtheorem{remark}[lemma]{Remark}
\newtheorem{definition}[lemma]{Definition}
\newtheorem{corollary}[lemma]{Corollary}

\begin{document}

\newcommand{\Lfr}{\mathfrak {L}}
\newcommand{\Afr}{\mathfrak {A}}
\newcommand{\Cpx}{\mathbf {C}}
\newcommand{\Fb}{\mathbf {F}}
\newcommand{\restrict}{\upharpoonright}
\newcommand{\svect}{\overset\rightarrow s}
\newcommand{\tvect}{\overset\rightarrow t}
\newcommand{\eps}{\epsilon}
\newcommand{\Fc}{\mathcal{F}}
\newcommand{\Mcal}{\mathcal{M}} 
\newcommand{\Nc}{\mathcal{N}}
\newcommand{\oup}{^{\mathrm o}}
\newcommand{\xh}{\hat {x}}
\newcommand{\yh}{\hat {y}}
\newcommand{\figf}{(\displaystyle\ast _{1}^{p}Q)\otimes M_{p}(\mathbb{C})}
\newcommand{\jconj}{\overline{\delta}}
\newcommand{\gruppo}{\mathbb{Z}_{p}}
\newcommand{\luna}{\mathbb{Z}_{p^2}}
\newcommand{\freeprodN}{\displaystyle \ast _{1}^{p}N}
\newcommand{\freeprodQ}{\displaystyle \ast _{1}^{p}Q}
\newcommand{\ifreeprod}{\left (\freeprodQ\right )\ast\mathfrak{L}\left (\Fb  _{1-\frac{1}{p}}\right )}
\newcommand{\pippo}{\left(\freeprodQ\ast \mathfrak{L}\left (\Fb _{1-\frac{1}{p}}\right )\right )\otimes R}
\newcommand{\ipippo}{\left (\ifreeprod\right )\otimes R}
\newcommand{\tenpro}{\left (\ifreeprod\right )\otimes R_{0}}
\newcommand{\pten}{\left (\tenpro\right )}
\newcommand{\crossprod}{\pten\rtimes_{\gamma}\gruppo}
\newcommand{\prodotto}{\left ( \ipippo\right )\rtimes_{\gamma}\gruppo} 
\newcommand{\obstruction}{e^{\frac{2\pi i}{p}}} 
\newcommand{\obstruconj}{e^{-\frac{2\pi i}{p}}} 
\newcommand{\jonesinv}{e^{\frac{2\pi i}{p^2}}}
\newcommand{\jonesconj}{e^{-\frac{2\pi i}{p^2}}}
\newcommand{\basic}{R_{-1}\rtimes _{\theta}\luna}
\newcommand{\dualgr}{\widehat{\gruppo}}
\newcommand{\M}{{\mathcal M}}
\newcommand{\dualcr}{\M\rtimes_{\widehat{\gamma}}\gruppo}
\newcommand{\freeproductBA}{\left (Q\otimes B\right )*_{B} A}
\newcommand{\freeproductDC}{\left (Q\otimes D\right )*_{D} C}
\newcommand{\freeproductRR}{\left (Q\otimes R_{0}\right )*_{R_{0}} R_{1}}
\newcommand{\freeprodRRo}{\left (Q\otimes R_{-1}\right )*_{R_{-1}} R_{0}}
\newcommand{\CD}{\mathcal{C}}
\newcommand{\tM}{\widetilde{M}}
\newcommand{\NQ}{N^{\CoD}(Q)}
\newcommand{\MQ}{M^{\CoD}(Q)}
\newcommand{\A}{{\mathcal A}}
\newcommand{\B}{{\mathcal B}}
\newcommand{\freeprolib}{(\libero\otimes N)*_{N}M}
\newcommand{\normalizer}{{\mathcal N}_{R_{0}}(R_{-1})}
\newcommand{\crossprodRo}{R_{0}\rtimes _{s}\gruppo}
\newcommand{\crossprodR}{R_{-1}\rtimes _{s}\gruppo}
\newcommand{\Ad}{\operatorname{Ad\,}}
\newcommand{\ad}{\operatorname{Ad}}
\newcommand{\comsquare}{\begin{array}{lcl}
 \vspace{.01cm}\hspace{5cm} B   &\vspace{.01cm}\subset &\vspace{.01cm}A \\
  \vspace{.01cm}\hspace{5cm}\cup &                      &\vspace{.01cm}\cup\\
   \hspace{5cm}  D                & \subset              & C      
                            \end{array}}

\newcommand{\qvd}{$\hfill\blacksquare$}
%
%PER AVERE LA STESSA SPAZIATURA
%
\newcommand{\inner}{\operatorname{Int}(\crossprodRo,\,\crossprodR)}
\newcommand{\innerM}{\operatorname{Int}(M)}
\newcommand{\innerN}{\operatorname{Int}(N)}
\newcommand{\innerP}{\operatorname{Int}(P)}
\newcommand{\innerR}{\operatorname{Int}(R)}
\newcommand{\innerA}{\operatorname{Int}(\A)}
\newcommand{\innerB}{\operatorname{Int}(\B)}
\newcommand{\innerfac}{\operatorname{Int}\left (\freefactor\right )}
\newcommand{\innerpten}{\operatorname{Int}\pten}
\newcommand{\CtM}{\operatorname{Ct}(M)}
\newcommand{\CtN}{\operatorname{Ct}(N)}
\newcommand{\CtP}{\operatorname{Ct}(P)}
\newcommand{\CtA}{\operatorname{Ct}(\A)}
\newcommand{\CtB}{\operatorname{Ct}(\B)}
\newcommand{\CtR}{\operatorname{Ct}(R_{0})}
\newcommand{\Ctpten}{\operatorname{Ct}\pten}
\newcommand{\Ctpi}{\operatorname{Ct}\left (\pippo\right )}
\newcommand{\autM}{\operatorname{Aut}(M)} 
\newcommand{\autR}{\operatorname{Aut}(R)}
\newcommand{\autpten}{\operatorname{Aut}\pten}
\newcommand{\autfac}{\operatorname{Aut}\left (\freefactor\right )}
\newcommand{\autN}{\operatorname{Aut}(N)} 
\newcommand{\autA}{\operatorname{Aut}(\A)}
\newcommand{\autB}{\operatorname{Aut}(\B)}
\newcommand{\outM}{\operatorname{Out}(M)}
\newcommand{\outN}{\operatorname{Out}(N)}
\newcommand{\outA}{\operatorname{Out}(\A)}
\newcommand{\outB}{\operatorname{Out}(\B)}
\newcommand{\outR}{\operatorname{Out}(R)}
\newcommand{\chiusoM}{\overline{\operatorname {Int}(M)}}
\newcommand{\chiusoN}{\overline{\operatorname{Int}(N)}}
\newcommand{\chiusoP}{\overline{\operatorname{Int}(P)}}
\newcommand{\chiusoR}{\overline{\operatorname{Int}(R)}}
\newcommand{\chiusofac}{\overline{\operatorname{Int}\left (\freefactor\right )}}
\newcommand{\chiusopten}{\overline{\operatorname{Int}\pten}}
\newcommand{\chiusoA}{\overline{\operatorname{Int}(\A)}}
\newcommand{\chiusoB}{\overline{\operatorname{Int}(\B)}}
\newcommand{\innerMN}{\operatorname{Int}(M,\, N)}
\newcommand{\innerRoR}{\operatorname{Int}(R_{0},R_{-1})}
\newcommand{\CtMN}{\operatorname{Ct}(M,\, N)}
\newcommand{\CtRoR}{\operatorname{Ct}(R_{0},\, R_{-1})}
\newcommand{\CtcrosRoR}{\operatorname{Ct}(\crossprodRo,\,\crossprodR)}
\newcommand{\autMN}{\operatorname{Aut}(M,\, N)} 
\newcommand{\autRoR}{\operatorname{Aut}(R_{0},R_{-1})}
\newcommand{\chiusoMN}{\overline{\operatorname {Int}(M,N)}}
\newcommand{\chiusoRoR}{\overline{\operatorname{Int}(R_{0},R_{-1})}}
\newcommand{\chiulibero}{\overline{\operatorname{Int}\left (\freeprolib\right )}}
\newcommand{\chiuso}{\overline{\operatorname{Int}(\crossprodRo,\, \crossprodR)}}

\pagestyle{myheadings}

\title[Existence of Non-Outer Conjugate Actions]{Non-Outer Conjugate $\luna$-Actions on Free Product Factors}
\author{Kenneth Dykema$^\dagger$ and Maria Grazia Viola}
\thanks{$\dagger$ Research partially supported by NSF grant DMS-0300336.}
\address{\hskip-\parindent Department of Mathematics \\Texas A\&M University \\ College Station TX 77843-3368, USA} 
\email{kdykema@math.tamu.edu}
\email{viola@math.tamu.edu}
\date{23 November 2004}

\begin{abstract}
We show that for any prime $p$ and for any II$_1$--factor $N$ 
there exist two $\mathbb{Z} _{p^2}$-actions on
the free product factor $*_1^pN$ that have the same outer invariant but are not outer conjugate.
Therefore, the outer invariant is not a complete invariant for outer conjugacy.
\end{abstract}

\maketitle

\section{Introduction}

Two automorphisms $\alpha$ and $\beta$ on a von Neumann algebra $M$ are said to be outer conjugate if there exists an automorphism $\sigma\in\autM$ and a unitary $W\in M$ such that $\sigma\beta\sigma ^{-1}=\Ad W\circ\alpha$, where $\Ad W$ acts on $M$ as $\Ad W (x)=WxW^{*}$.

In an exceptional paper \cite{Connes2} A. Connes classified periodic automorphisms on the hyperfinite $II_{1}$ factor $R$ up to outer conjugacy. He showed 
that the pair $(p_{0},\lambda)$, where $p_{0}$ is the outer period of the automorphism $\alpha\in\operatorname{Aut}(R)$, i.e., the smallest non-negative integer such that $\alpha ^{n}$ is an inner automorphism, and $\lambda$ is the obstruction to lifting (as defined in Section 2), is a complete outer conjugacy invariant for automorphisms. The pair $(p_{0},\lambda)$ is usually referred to as the outer invariant of the automorphism $\alpha$. 

The question addressed in this paper is whether similar results hold for other type of II$_{1}$ factors, in particular for free product factors. In fact, we show that the situation for free product factors is completely different. Given any prime $p$ and a II$_1$ factor $N$, we show that there exist two actions of the cyclic group $\luna$ on $\freeprodN$ which have the same outer period and obstruction to lifting, but are not outer conjugate. This work is a generalization of the work of Florin R\u{a}dulescu (see \cite{Radu2}) for $\mathbb{Z} _{2}$-kernels on an interpolated free group factor.

The two $\luna$-actions we consider are realized starting with
two very different $II_{1}$ factors $\A$ and $\M$ with , each of which is obtained from a $II_{1}$ subfactor construction.
The first factor $\A$ is an amalgamated free product and is generated by an extremal commuting square of finite dimensional algebras. We prove in Theorem \ref{theorem4.8} the existence of a pair $A_{1}\subseteq B_{1}$ of II$_1$ factors such that $\A$ is the enveloping algebra of this pair. Here, $B_{1}$ is the crossed product of $A_{1}$ by a $\luna$--action $\theta_{1}$ with outer invariant $(p, \obstruction)$, and $A_{1}$ is a free product factor.

The second factor is the crossed product $\M =\prodotto$, which is an example of a $II_{1}$ factor non-antiisomorphic to itself. $\M$ is the enveloping algebra of a pair $A_{2}\subseteq B_{2}$ with properties similar to those of the pair $A_{1}\subseteq B_{1}$ considered above.
 Here $Q$ is chosen so that $Q_{t}=N\otimes M_{p}(\mathbb{C})$ for $t=\sqrt{\frac{p+1}{p^3}}$.
 In particular $B_{2}$ is the crossed product of $A_{2}$ by a $\luna$--action  $\theta_{2}$ with outer invariant $(p, \obstruction)$, and $\displaystyle A_{1}\cong A_{2}$. 
 
In addition, letting $q$ be a projection of trace $t=\sqrt{\frac{p+1}{p^3}}$, we get $qA_{1}q\cong\freeprodN$. Therefore, perturbing $\theta _{i}$, i=1,2, by an inner automorphism $\phi _{i}$, $i=1,2$, so that $\phi _{i}\circ\theta _{i}$ i=1,2, leave $q$ invariant, we obtain two $\luna$-actions $\widetilde{\theta _{i}}=\phi _{i}\circ\theta_{i}|_{qA_{1}q}$ on the free product factor $\freeprodN$ that have also outer invariant $(p, \obstruction)$. 

Our subfactor constructions thus give two $\luna$-actions on $\freeprodN$ with the same outer invariant. To show that these two actions are not outer conjugate we compute the Connes $\chi$ invariant of the two factors $\M$ and $\A$. This technique was already used by R\u{a}dulescu in \cite{Radu2} to show that the $\mathbb{Z} _{4}$-actions he constructed were not outer conjugate. The Connes invariant, introduced by Connes in \cite{Connes0}, is a certain abelian subgroup of the group of outer automorphisms of a $II_{1}$ factor, as recalled  in Section 2. Although for both factors the Connes invariant is algebraically isomorphic to $\luna$, it is still
able to distinguish the factors if we view the invariant
as a subgroup of the group of outer automorphisms. In particular, we are interested in the position of the unique subgroup of order $p$ of $\chi (\M)\cong\chi(\A)\cong\luna$. We 
thus consider the crossed product of each of the factors $\M$ and $\A$ by the $\gruppo$-action corresponding to this unique subgroup of the $\chi$ invariant. As shown in Section 5 the associated dual actions can be decomposed, modulo inner automorphisms, into an approximately inner automorphism and a centrally trivial automorphism. Using this decomposition we prove that the two dual $\gruppo$--actions are not outer conjugate, by showing that the two factors $\M$ and $\A$ are not
not isomorphic (Theorem 6.5). It follows that the two $\luna$-actions on $\freeprodN$ are not outer conjugate.

\section{Definitions}

Let $M$ be a $II_{1}$ factor with separable predual, endowed with a faithful trace $\tau$. Then $M$ inherits an $L^{2}$-norm from the inclusion $M\subset L^{2}(M)$ given by $\|x\|_{2}=\tau(x^{*}x)^{\frac{1}{2}}$, for all $x\in M$. Denote by $\autM$ the group of automorphisms of M endowed with the pointwise weak-topology for which a sequence of automorphisms $\alpha _{n}$ converges to $\alpha$ if and only if\, $\|\alpha _{n}(x)-\alpha (x)\|_{2}\rightarrow 0$ for all $x\in M$.

When studying the automorphisms of a factor there are two normal subgroups of $\autM$ which are of particular interest, the group of inner automorphisms and the group of centrally trivial automorphisms. These are also the groups involved in the definition of the Connes invariant $\chi (M)$, as we will briefly describe. The inner automorphisms are the automorphisms of the form $\ad _M (u)$ for some unitary $u$ in $M$. If there is no ambiguity in determining the factor M on which the inner automorphism acts, we will simply use the notation $\Ad u$. Let $\innerM$ be the normal subgroup of $\autM$ formed by all inner automorphisms, and denote by $\chiusoM$ its closure in the pointwise weak-topology. 

\begin{definition}
\label{defn2.1}
A bounded sequence $(x_{n})_{n\in\mathbb{N}}$ in $M$ is called {\em central} if \linebreak $\displaystyle{\lim_{n \rightarrow \infty}{\| x_{n}y-yx_{n}\|_{2}=0}}$, for all $y\in M$. 
\end{definition}

Also, a central sequence $(x_{n})_{n\in\mathbb{N}}$ in $M$ is trivial if there exists a sequence of complex numbers $(\lambda _{n})_{n\in\mathbb{N}}$ such that $\displaystyle{\lim _{n\rightarrow\infty}\|x_{n}-\lambda _{n}1\|_{2}=0}$.

Next we define the centrally trivial automorphisms of $M$, which form the other group involved in the definition of the Connes invariant.   

\begin{definition} 
\label{defn2.2}
An automorphism $\alpha\in\autM$ is said to be {\em centrally trivial} if $\displaystyle{\lim_{n \rightarrow \infty}{\| \alpha (x_n)-x_n \|_{2}}=0}$ for any central sequence $(x_{n})_{n\in\mathbb{N}}$ in $M$.  
\end{definition}

\noindent Let $\CtM$ be the group of centrally trivial automorphisms
and $\displaystyle{\outM=\frac{\autM}{\innerM}}$ be the group of outer automorphisms of $M$,
with
\[
\epsilon:\autM\longrightarrow\outM
\]
the quotient map.
We can now define the Connes invariant $\chi (M)$, which was introduced by Connes in \cite{Connes0}.
 
\begin{definition}
\label{defn2.3}
   Let M be a $II_{1}$ factor with separable predual. The {\em Connes invariant} $\chi (M)$ is the abelian group 
   \[\chi(M) =\frac{\CtM\cap\chiusoM}{\innerM}\subset\outM.\]
\end{definition}

In a similar fashion, given an inclusion $N\subseteq M$ of $II_{1}$ factors with separable predual, one can define the relative Connes invariant in the following way. Let $\autMN$ be the group of automorphisms of $M$ leaving $N$ invariant and let $\innerMN$ be the subgroup of $\autMN$ formed by all inner automorphisms of $M$ implemented by unitaries in $N$. Denote by $\CtMN$ the set of automorphisms of $M$ leaving $N$ invariant and acting asymptotically trivially on the central sequences for $M$ which are contained in 
$N$. The relative Connes invariant for the inclusion $N\subseteq M$ was introduced by Y. Kawahigashi in \cite{Kawi} as a generalization of the Connes $\chi$ invariant and is defined as  
\begin{equation*}
\chi (M,N)=\frac{\CtMN\cap\chiusoMN}{\innerMN}.
\end{equation*} 

The most useful tool for studying the Connes invariant or its relative version is the central sequence algebra. Let $\omega$ be a free ultrafilter over $\mathbb{N}$ and denote by $\tau$ the trace on the $II_1$ factor $M$. As usual we assume that $M$ has separable predual. Let $\ell ^{\infty}(\mathbb{N},M)$ be the algebra of bounded sequences in $M$. Define $\mathfrak{I}_{\omega}$ as the ideal in $\ell ^{\infty}(\mathbb{N}, M)$ formed by all sequences $(x_{n})_{n\in\mathbb{N}}$ such that $\displaystyle{\lim _{n\rightarrow\omega}\|x_{n}\|_{2}=0}$ (see \cite{EvK} for more details on ultrafilters). Set $\displaystyle{M^{\omega}=\frac{\ell ^{\infty}(\mathbb{N},M)}{\mathfrak{I}_{\omega}}}$. Then $M$ can be embedded into $M^{\omega}$ as the set of constant sequences. 

Let $C_{\omega}$ denote the subalgebra of $\ell ^{\infty}(\mathbb{N},M)$ consisting of sequences with the property that $\displaystyle{\lim _{n\rightarrow\omega}\|x_{n}y-yx_{n}\|_{2}=0}$ for all $y\in M$. We will refer to these sequences as $\omega$-central.

\begin{definition}
\label{defn2.4}
For a free ultrafilter $\omega$ over $\mathbb{N}$, the {\em $\omega$--central sequence algebra} of $M$ is 
\[M_{\omega}=\frac{C_{\omega}}{\mathfrak{I}_{\omega}\cap C_{\omega}}.\]
\end{definition}

\noindent Set $\displaystyle\tau ^{\omega}([(x_{n})_{n}])=\lim_{n\rightarrow\omega}\tau(x_{n})$, where $[(x_{n})_{n}]$ denotes the coset $(x_{n})_{n}+\mathfrak{I}_{\omega}$ in $M^{\omega}$ and $(x_{n})_{n}\in\ell ^{\infty}(\mathbb{N},M)$. Then $\tau^{\omega}$ defines a faithful trace on $M^{\omega}$, which induces a faithful trace on $M_{\omega}$. Moreover, $M_{\omega}=M^{\omega}\cap M^{\prime}$. 

Observe that every automorphism $\alpha\in\autM$ induces an automorphism $\alpha^{\omega}$ on $M^{\omega}$ by 
\[\alpha^{\omega}([(x_{n})_{n}])=[(\alpha(x_{n}))_{n}],\]
for all $[(x_{n})]\in M^{\omega}$. By restricting this automorphism to classes of $\omega$--central sequences we obtain an automorphism $\alpha _{\omega}$ of $M_{\omega}$ for every $\alpha\in\autM$. The following remark, due to Connes \cite{Connes1}, gives a characterization of the centrally trivial automorphisms.

\begin{remark}
\label{remark2.5}
Let $N\subseteq M$ an inclusion of $II_{1}$ factors with separable predual. An automorphism $\alpha\in\autMN$ belongs to $\CtMN$ if and only if 
\[\alpha _{\omega}|_{N^{\omega}\cap M^{\prime}}=Id|_{N^{\omega}\cap M^{\prime}}\]
\end{remark}

The Connes classification of the periodic automorphisms on the hyperfinite $II_{1}$ factor, up to outer conjugacy, is based on the outer invariant of an automorphism $\alpha$ of $M$, which is defined as follows. 

\begin{definition}
\label{defn2.6}
Given a $II_{1}$ factor $M$, the {\em outer period} $p_{0}$ of $\alpha\in\autM$ is the smallest non-negative integer such that 
$\alpha ^{p_{0}}=\ad _{M}(U)$, for some unitary $U\in M$. If no power of $\alpha$ is an inner automorphism then we set $p_{0}=0$. Note that if $U$ is
a unitary that realizes $\alpha ^{p_{0}}$ as an inner automorphism, then it is easy to check that $\alpha (U)=\lambda\,U$, for some $\lambda\in\mathbb{C}$ with $\lambda^{p_{0}}=1$ (apply $\alpha^{p_{0}+1}=\alpha ^{p_{0}}\alpha=\alpha\alpha^{p_{0}}$ to $x$ in $M$ and use the fact that $M$ is a factor). We call $\lambda$ the {\em obstruction to lifting} of $\alpha$ and the pair $(p_{0}, \lambda)$ the {\em outer invariant} of $\alpha$. We sometimes refer to a periodic automorphism of $M$ with outer period $p$ as a $\gruppo$-kernel. 
\end{definition}

The terminology obstruction to lifting is motivated by the fact that $\lambda=1$ if and only if the homomorphism $\phi :\mathbb{Z}_{p_{0}}\longrightarrow\outM$ with $\phi (1)=\alpha$ can be lifted to an homomorphism $\Phi :\mathbb{Z}_{p_{0}}\longrightarrow\autM$, so that the diagram 

\begin{picture}(100,120)(0,0)
\put (200,10){$\mathbb{Z}_{p_{0}}$\vector(1,0){37}$\outM$}
\put (230,14){$\phi$}
\put (215,20){\vector(3,4){50}}
\put (225,53){$\Phi$}
\put (267,85){$\autM$}
\put (280,80){\vector(0,-1){60}}
\put (282,50){$\varepsilon$}
\end{picture}

commutes.

\section{Commuting Squares}

In this section we review some basic elements about commuting squares as they appear in \cite{Popa1}, \cite{Popa3} (see also \cite{GHJ}). Given a set $S\subseteq B(H)$ denote by $\mbox{Span(S)}$ the linear span of $S$, and by $\overline{\mbox{Span}(S)}^{\text{wo}}$ its closure in the weak operator topology.

\begin{definition}[Popa]
\label{defn3.1}
A diagram  
\[\comsquare\]
of finite von Neumann algebras with a finite faithful normal trace $\tau$ on 
$A$ is a {\em commuting square} if the diagram 
\[\begin{array}{lcl}
     \vspace{.01cm}\hspace{5cm} B & \vspace{.01cm}\stackrel{E_{B}}{\longleftarrow} & 
     \vspace{.01cm} A    \\
     \vspace{.01cm}\hspace{5cm} \uparrow\it{i} &    & \vspace{.01cm}\uparrow\it{i} \\
    \hspace{5cm}  D      & \stackrel{E_{D}}{\longleftarrow} & C      
\end{array}\]
commutes, where $i$ denotes the inclusion map, and $E_{D}$, $E_{B}$ are the conditional expectations onto $D$ and $B$ respectively.
\end{definition}

The theory of von Neumann algebras and subfactors provides numerous examples of commuting squares. In fact, if $N\subset M\subset M_{1}\subset M_{2}\subset\cdots$\; is the Jones tower for the inclusion $N\subset M$, then 
for all $i\ge1$ the relative commutants
\[\begin{array}{ccc}
   \vspace{.01cm}\hspace{5cm} M^{\prime}\cap M_{i+1} & \vspace{.01cm}\subset & \vspace{.01cm} M_{1}^{\prime}\cap M_{i+1}    \\
   \vspace{.01cm}\hspace{5cm} \cup&    & \vspace{.01cm}\cup  \\
    \hspace{5cm}  M^{\prime}\cap M_{i}  & \subset & M_{1}^{\prime}\cap M_{i}      
\end{array}\] 
form a commuting square.

We are interested in commuting squares that have some additional properties, which appear naturally when studying commuting squares coming from a finite depth inclusion of subfactors.

\begin{definition}
\label{defn3.2}

    Let 
    \[\comsquare\]
    be a commuting square of finite von Neumann algebras 
    which are weakly separable. Denote by $\tau$ the normalized, faithful trace defined on 
    $A$. We say that the commuting square is {\em extremal} if it satisfies the following conditions:
    \begin{description}
    \item[(i)] The inclusion $C\subset A$ has a $\lambda$-Markov trace \cite{GHJ},
      (which means that there exists a trace $\widetilde{\tau}$ on 
      the basic construction \[ \langle A, f\rangle =\overline{Span(AfA)}^{\text{wo}}\] 
      for this inclusion that extends the trace $\tau$ on $A$, and has the property that  
      $\widetilde{\tau} (fa)=\lambda\tau (a)$ for all $a\in A$)
      such that the 
      basic construction for $D\subseteq B$, with respect to the trace 
      $\widetilde{\tau}|_{B}$, coincides with $\langle B,f\rangle$
   \item[(ii)] $C$ generates $A$ as a left (or right) B-module, i.e.,
         \[ A=\overline{Span (BC)}^{\text{wo}}=\overline{Span (CB)}^{\text{wo}} \]
   \item[(iii)] The centers of the pair of algebras $A,B$
      and $C,D$, respectively, have trivial intersection  
      \[{\mathcal Z}(A)\cap {\mathcal Z}(B)=\mathbb{C}\, 1\quad\text{and}\quad{\mathcal Z}(C)\cap 
      {\mathcal Z}(D)=\mathbb{C}\, 1.\]             
     This is equivalent to saying that the inclusions $B\subset A$ and 
$D\subset C$ have connected Bratteli diagrams.
\end{description}
\end{definition}

Note that if $A$ is a finite dimensional algebra the condition (ii) is unnecessary because equivalent to (i). There is a standard way, the basic construction, to produce a new extremal commuting square from a given one. 

\begin{remark}
\label{remark3.3}
Given an extremal commuting square of finite dimensional algebras 
\[\comsquare \]
construct the diagram 
\[\begin{array}{ccc}
     \vspace{.01cm}\hspace{5cm} \langle B, f\rangle & \vspace{.01cm}\subset 
     & \vspace{.01cm} \langle A, f\rangle \\
     \vspace{.01cm}\hspace{5cm} \cup&    & \vspace{.01cm}\cup  \\
     \hspace{5cm}  B  & \subset & A     
\end{array}\] 
This is also an extremal commuting square. 
\end{remark}

\noindent As mentioned before, finite depth inclusions of factors provide many examples of extremal commuting squares. In fact, it is a result of A. Ocneanu \cite{Ocneanu2} that if $N\subseteq M$ is an arbitrary inclusion of $II_{1}$ hyperfinite factors of finite depth then there exists an extremal commuting square of finite dimensional algebras which by iteration of the basic construction produces the inclusion $N\subseteq M$. 

\section{The construction of the two $\luna$-actions}

Fix a prime $p$ and let $N$ be a II$_1$ factor.
In this section we construct two $\luna$-actions on $\freeprodN$ which have the same outer invariant. We will prove later (Section 6) that these actions are not outer conjugate. 

Denote by $R$ the hyperfinite $II_{1}$ factor. The first one of our actions is obtained from
an inclusion $A_{1}\subset B_{1}$ of $II_{1}$ factors which by iteration of the basic construction produces
the $II_{1}$ factor $(Q\otimes R)*_{R}(R\rtimes _{\sigma}\luna)$, where $Q$ is a $II_{1}$ factor and $\sigma$ is a $\luna$-action on $R$. In Theorem 4.4 we show that $B_{1}$ is the crossed product of $A_{1}$ by a $\luna$-action with outer invariant $(p,\obstruction)$. We will denote the action we define on $\freeprodN$ by $\widetilde{\theta_{1}}$. 
Analogously, the second action is obtained from an inclusion $A_{2}\subset B_{2}$, which by iteration of the basic construction produces the II$_{1}$ factor $\crossprod$ and has similar properties than the inclusion $A_{1}\subset B_{1}$.

 To construct the two $\luna$-actions and to prove that they are not outer conjugate it is  useful to have an explicit model for the hyperfinite $II_{1}$ factor $R$ and for the inclusion $R\subset R\rtimes _{\sigma}\luna$. Such model was introduced by R\u{a}dulescu in \cite{Radu2} for the case $p=2$. Since the proof for the general case is based on the same argument given by R\u{a}dulescu in his work, we refer to Lemma 5 in \cite{Radu2} for a proof of the following statement.
   
\begin{lemma}
\label{lemma4.1}
The hyperfinite II$_1$--factor is generated by unitaries $(U_k)_{k\in\mathbb{Z}}$ and a unitary $g$,
where each $U_k$ has spectrum consisting of the $p^2$-th roots of unity, on which the trace of spectral measure
is equidistributed, and where $g$ has spectrum consisting of the $p$-th roots of unity with equidistributed
trace of spectral measure, where the following relations hold:
\begin{itemize}
\item [(i)]  $(U_{k})^{p^2} =1\text{ for all }k\in\mathbb{Z}$,
\item [(ii)]
           $U_{k}gU_{k}^* =\obstruconj g\text{ if }k=0,-1,
           \text{ while } U_{k}gU_{k}^{*}=g,\text{ if }k\in\mathbb{Z}\backslash\{0,-1\}$,
\item [(iii)] $U_{k} U_{k+1}U_{k}^{*}=\jonesinv U_{k+1},\text{ for } 
              k\in\mathbb{Z}$, 
\item [(iv)]
            $U_{i}U_{j}=U_{j}U_{i}\text{ if }\quad|i-j|\geq 2$,
\end{itemize}   
and where the trace of every nontrivial word in the $U_k$ and $g$ (in the obvious sense, subject to the above relations)
is zero.

Moreover,
$$
R_{-1}=\{gU_{0}^{p}, U_{1}, U_{2},\hdots\}^{\prime\prime}\subseteq\{g,U_{0}, U_{1}, U_{2},\hdots\}^{\prime\prime}=R_{0}.
$$
defines an inclusion of type $\Pi _1$ factors of index $p^2$ for which the relative commutant $R_{-1}'\cap R_{0}$ is equal to $\{g\}''$.

If we set $\theta=Ad_{R_{-1}}(U_0)$, then $\theta$ is a $\luna$-action of $R_{-1}$ and $R_0$ is equal to the crossed product $\basic$. Moreover, the outer invariant of 
$\theta$ is $(p,\obstruconj)$. 

Also, the tower of factors in the iterated basic construction for the inclusion  $R_{-1}\subset R_0$ is given by the family
   $$
   R_{k} =\{g, U_{-k}, U_{-k+1}, U_{-k+2}, \hdots\}^{\prime\prime}\,\text{ for }k\geq 1, 
   $$
while the Jones downward tunnel for $R_{-1}\subset R_0$ is given by
$$
 R_{-k} =\{gU_{0} ^{p}U_{1} ^{p}\cdots U_{k-1} ^{p}, U_{k}, U_{k+1},U_{k+2},\hdots\}^{\prime\prime}\,\text{ for } k\geq 1, 
   $$
Moreover, the k-th term in the sequence of relative commutants is
$$
R_{k}\cap R_{-1}^{\prime}=\{g, U_{-k}, U_{-k+1},\hdots ,U_{-1}\}^{\prime\prime},\,\text{ for }k\geq 1.
$$
\end{lemma}

Since $R_{0}=\{g,U_{0}, U_{1}, U_{2},\hdots\}^{\prime\prime}$ is obviously generated by an increasing family of finite dimensional subalgebras, it provides a model for the hyperfinite $II_{1}$ factor. Note also that it is quite simple to produce a family of unitaries satisfying the conditions of Lemma 4.1. In fact, if $\beta$ is a $\luna$-action on the hyperfinite $II_{1}$ factor with outer invariant $(p, \obstruction)$, then the unitaries implementing the crossed products in the Jones
tower and tunnel for the inclusion $R\subset R\rtimes_{\beta}\luna$ verify the desired properties. The unitary $g$ of order $p$ is chosen instead from the elements of the first relative commutant. 

 In the next lemma we show that the inclusion $R_{-1}\subset\basic =R_{0}$ can be obtained by iteration of the basic construction from an extremal commuting square (see comment at the end of Section 3). 

\begin{lemma}
\label{lemma4.2}
  Set $R_{0}=\{g,U_{0}, U_{1}, U_{2},\hdots\}^{\prime\prime}\subset R_{1}=\{g,U_{-1}, U_{0}, U_{1}, U_{2},\hdots\}^{\prime\prime}$, and $R_{-k} =\{gU_{0} ^{p}U_{1} ^{p}\cdots U_{k-1} ^{p}, U_{k}, U_{k+1},U_{k+2},\hdots\}^{\prime\prime}$ for $k\geq 1$, as in the previous lemma. 
Then
\[
\begin{array}{rcl}
\{g, U_{0}\}^{\prime\prime}=B & \subset & A=\{g, U_{-1},U_{0}\}^{\prime\prime} \\
\cup & & \cup \\
\{g\}^{\prime\prime}=D & \subset & C=\{g, U_{-1}\}^{\prime\prime}.
\end{array}
\]
is an extremal commuting square and generates the inclusion $R_{0}\subset R_{1}$
by iteration of the Jones basic construction for $D\subset B$ and $C\subset A$.
The successive steps in the Jones tower for the inclusions $D\subset B$ and $C\subset A$ are given by the algebras $A_{0} ^{(k)}$ and $A_{1} ^{(k)}$ respectively, 
   where $A_{0} ^{(k)}$ and $A_{1} ^{(k)}$ are defined by
\[
 A_{0} ^{(k)}=\{g, U_{0}, U_{1}, \hdots,U_{k}\}^{\prime\prime}\subset
 A_{1} ^{(k)}=\{g, U_{-1}, U_{0}, U_{1},\hdots, U_{k}\}^{\prime\prime}
\]
for $k\ge0$ and
\[
A_{0} ^{(-1)}=D=\{g\}''\subset A_{1} ^{(-1)}=C=\{g,U_{-1}\}''.
\]
Moreover, $A_{0} ^{(k-2)}=R_{-k}'\cap R_{0}$ and  $A_{1} ^{(k-2)}=R_{-k}'\cap R_{1}$ for all $k\geq 1$
\end{lemma}

\begin{proof}
The diagram  
  \[\begin{array}{ccc}
     \vspace{.01cm}\hspace{5cm} \{g, U_{0}\}^{\prime\prime} & \vspace{.01cm}\subset 
     & \vspace{.01cm} \{g,U_{-1}, U_{0}\}^{\prime\prime} \\
     \vspace{.01cm}\hspace{5cm} \cup&    & \vspace{.01cm}\cup  \\
     \hspace{5cm}  \{g\}^{\prime\prime}  & \subset & \{g,U_{-1}\}^{\prime\prime}.     
  \end{array}\]            
is an extremal commuting square because of the properties satisfied by the unitaries $U_{k}$ and $g$
(i.e.\ (ii)-(iv) in the previous lemma) and the definition of the trace on $R_{0}$.

Next we want to show that $A_{0} ^{(k-2)}=R_{-k}'\cap R_{0}\subset R_{-k}'\cap R_{1}=A_{1} ^{(k-2)}$. Set 
\begin{eqnarray*}
& & A_{-j}^{(j-1)} = \{gU_{0}^{p}U_{1}^{p}\cdots U_{j-1}^{p}\}^{\prime\prime}\\
& & A_{-j}^{(k)}=\{gU_{0}^{p}U_{1}^{p}\cdots U_{j-1}^{p}, U_{j}, U_{j+1}, U_{j+2},\hdots U_{k}\}^{\prime\prime}\text{ for }k\geq j\text{ and }j\geq 1.
\end{eqnarray*}

By the previous lemma $U_{k}$ is the unitary implementing the crossed product in the $k$-th step of the Jones tunnel for the inclusion $R_{0}\subset R_{0}\rtimes _{Ad(U_{-1)}}\luna =R_{1}$, so iterating the basic construction horizontally and constructing the Jones tunnel vertically, we obtain the following diagram 

\[\begin{array}{llllllclllll}
 A_{1} ^{(-1)}\subset & A_{1} ^{(0)}\subset & A_{1} ^{(1)}\subset & A_{1} ^{(2)}\subset &\hdots & A_{1} ^{(k-1)}\subset & A_{1} ^{(k)}\subset & A_{1} ^{(k+1)}\subset & \hdots \subset &R_{1}\\  
\cup & \cup & \cup & \cup &  & \cup & \cup & \cup & & \cup \\
A_{0} ^{(-1)}\subset & A_{0} ^{(0)}\subset & A_{0} ^{(1)}\subset & A_{0} ^{(2)}\subset &\hdots & A_{0} ^{(k-1)}\subset & A_{0} ^{(k)}\subset & A_{0} ^{(k+1)}\subset & \hdots \subset &R_{0}\\
  & \cup & \cup & \cup & \hdots & \cup & \cup & \cup & & \cup \\
 & A_{-1} ^{(0)}\subset & A_{-1} ^{(1)}\subset & A_{-1} ^{(2)}\subset &\hdots & A_{-1} ^{(k-1)}\subset & A_{-1} ^{(k)}\subset & A_{-1} ^{(k+1)}\subset & \hdots\subset &R_{-1}\\
 & & \cup & \cup & \hdots & \cup & \cup & \cup & &\cup \\
 & & A_{-2} ^{(1)}\subset & A_{-2} ^{(2)}\subset &\hdots & A_{-2} ^{(k-1)}\subset & A_{-2} ^{(k)}\subset & A_{-2} ^{(k+1)}\subset & \hdots\subset &R_{-2}\\
 & & & \cup & & \cup & \cup & \cup & &\cup \\
 & & &  \vdots & & \vdots & \vdots & \vdots & & \vdots \\
 & & & &  & \cup & \cup & \cup &  \,&  \cup \\
 & & & & & A_{-k} ^{(k-1)}\subset & A_{-k} ^{(k)}\subset & A_{-k} ^{(k+1)}\subset & \hdots \subset &R_{-k}\\ 
 & & & & & \cup & \cup & \cup &  &  \cup \\
 & & & & & \vdots & \vdots & \vdots & &\vdots 
\end{array}\]
where $A_{-i} ^{(j)}=A_{-i} ^{(j-1)}\rtimes_{Ad(U_{j})}\luna$ for all $i\geq 0$ and $j\geq i$, as well as for $i=1$ and $j\geq 0$.     

Obviously  $A_{0} ^{(k-2)}\subset R_{-k}'\cap R_{0}$, so we only need to show the other inclusion. Let $x\in R_{-k}'\cap R_{0}$. Since $\{A_{0} ^{(j)}\}_{j\geq -1}$ is the Jones tower of $R_{0}$, if we denote by $E_{A_{0}^{(j)}}^{R_{0}}$ the conditional expectation of $R_{0}$ onto $A_{0}^{(j)}$, then $x=\displaystyle \lim_{j\rightarrow\infty}x^{(j)}$ with $x^{(j)}=E_{A_{0}^{(j)}}^{R_{0}}(x)$.

Since $x\in R_{-k}^{\prime}$ we have that for any $y\in A_{-k}^{(j)}\subset A_{0}^{(j)}$  
\[yE_{A_{0}^{(j)}}^{R_{0}}(x)=E_{A_{0}^{(j)}}^{R_{0}}(yx)=E_{A_{0}^{(j)}}^{R_{0}}(xy)=E_{A_{0}^{(j)}}^{R_{0}}(x)y.\]
Therefore $x^{(j)}\in (A_{-k}^{(j)})'\cap A_{0}^{(j)}$ for every $j\geq k-1$.

Since $A_{0} ^{(j)}=A_{0} ^{(j-1)}\rtimes_{Ad(U_{j})}\luna$ for all $j\geq 0$, by iterating the crossed product construction, we can write any element $x^{(j)}\in (A_{-k}^{(j)})'\cap A_{0}^{(j)}$ as $x^{(j)}=\displaystyle\sum_{\overline {l}}a_{\overline{l}}W^{\overline{l}}$, where $\overline{l}=(l (k-1), l (k),\hdots, l (j))$ is an element of  $\{1,2,\hdots, p^2 -1\}^{j-k+2}$, $a_{\overline{l}}\in\{g, U_{0}, U_{1},\hdots ,U_{k-2}\}''=A_{0} ^{(k-2)}$ and $W^{\overline{l}}=U_{k-1}^{l (k-1)}U_{k}^{l (k)}U_{k+1}^{l (k+1)}\cdots U_{j}^{l (j)}$. Now we want to explore the restrictions given by the condition $x^{(j)}\in (A_{-k} ^{(j)})'$ for $j\geq k$. A simple computation shows that 
\begin{itemize}
\item[]\hspace{-1.1cm} $Ad(gU_{0} ^{p}U_{1} ^{p}\cdots U_{k-1} ^{p})(W^{\overline{l}})=e^{\frac{2\pi i(l(k)-l (k-1))}{p}}W^{\overline{l}}$  
\item[]\hspace{-1.1cm} $Ad(U_{m})(W^{\overline{l}})=e^{\frac{2\pi i(l(m+1)-l (m-1))}{p^2}}W^{\overline{l}}$ for all $k\leq m\leq j-1$  
\item[]\hspace{-1.1cm}$Ad(U_{j})(W^{\overline{l}})=e^{-\frac{2\pi i l (j-1)}{p^2}}W^{\overline{l}}$. 
\end{itemize}
Thus, $x^{(j)}\in {(A_{-k}^{(j)}})'\cap {A_{0}^{(j)}}$ if and only if
\begin{itemize}
\item []\hspace{-1.1cm}   i) $l (j-1)=0$
\item []\hspace{-1.1cm}  ii) $l (m+1)=l (m-1)$ for all $k\leq m\leq j-1$   
\item []\hspace{-1.1cm} iii) $l (k)\equiv l (k-1)\bmod p$.
\end{itemize}
Assume $j-k\equiv 0\bmod 2$. Then, the previous conditions imply that $x^{(j)}\in (A_{0} ^{(k-2)}\cup\{U_{k} ^{p}U_{k+2} ^{p}U_{k+4} ^{p}\cdots U_{j-2} ^{p}U_{j} ^{p}\})''$, with $[U_{k} ^{p}U_{k+2} ^{p}U_{k+4} ^{p}\cdots U_{j-2} ^{p}U_{j} ^{p}, A_{0} ^{(k-2)}]=0$. Since \linebreak $x^{(j)}=E_{A_{0}^{(j)}}^{R_{0}}(x)\in (A_{-k}^{(j)})'\cap A_{0}^{(j)}$ we also have that $x^{(j)}=E_ {(A_{-k} ^{(j)})'\cap A_{0}^{(j)}}(x^{(j+2)})$, where \linebreak $x^{(j)}=\displaystyle\sum_{l=0} ^{p-1}a_{l}(U_{k} ^{p}U_{k+2} ^{p}U_{k+4} ^{p}\cdots U_{j-2} ^{p}U_{j} ^{p}U_{j+2} ^{p})^{l}$ and $a_{l}\in A_{0} ^{(k-2)}$. Note that for every $l\neq 0$ $(U_{k} ^{p}U_{k+2} ^{p}U_{k+4} ^{p}\cdots U_{j-2} ^{p}U_{j} ^{p}U_{j+2} ^{p})^{l}$ is orthogonal to $(A_{-k}^{(j)})'\cap A_{0}^{(j)}$ with respect to $E_ {(A_{-k}^{(j)})'\cap A_{0}^{(j)}}$, thus
$$
x^{(j)}=E_ {(A_{-k}^{(j)})'\cap A_{0}^{(j)}}(x^{(j+2)}=a_{0}\in  A_{0} ^{(k-2)}.
$$ 
Similarly, if $j-k\equiv 1\bmod 2$ one can show that  $x^{(j)}\in (A_{0} ^{(k-2)}\cup\{U_{k-1} ^{p}U_{k+1} ^{p}\cdots U_{j} ^{p}\})^{\prime\prime}$ belongs to $A_{0} ^{(k-2)}$. Since $x=\displaystyle\lim_{j\rightarrow\infty}x^{(j)}$, it follows that $x\in  A_{0} ^{(k-2)}$, so $R_{-k}'\cap R_{0}\subset  A_{0} ^{(k-2)}$. An analogous argument shows that  $R_{-k}'\cap R_{1}=A_{1} ^{(k-2)}$ .
\end{proof}

Using the above commuting square we can now prove the existence of an inclusion of von Neumann algebras $A_{1}\subset B_{1}$ for the amalgamated free product $\A=(\freeproductRR)$, where $Q$ is a $II_{1}$ factor. Moreover, $B_{1}$ is the crossed product of $A_{1}$ by a $\luna$-action, and $A\cong (*_{1}^{p}Q)\otimes M_{p}(\mathbb{C})$. The proof uses the following lemma, which is a consequence of~\cite{Dy}.

\begin{lemma}
\label{lemma4.3}
Let $N\subset M$ be an irreducible inclusion of finite dimensional von Neumann algebras
with inclusion matrix $\Lambda$ and let $\tau$ be a faithful, tracial state on $M$ with trace vector
$\svect$.
Let
\[
\tvect=(t_j)_{j=1}^\ell=\Lambda^T\svect
\]
denote the trace vector for $\tau |_N$.
Thus,
\[
N=\bigoplus_{j=1}^\ell M_{n_j}(\Cpx)
\]
and $t_j$ is the result of $\tau$ applied to a minimal projection of the $j$th summand
of $N$.
Let $Q$ be any II$_1$--factor.
Then for the amalgamated free product of von Neumann algebras with respect to the
obvious trace--preserving conditional expectations onto $N$, we have
\begin{equation}\label{eq:QNM}
(Q\otimes N)*_NM\cong Q_{\frac1{t_1}}*\cdots*Q_{\frac1{t_\ell}}*L(\Fb_r),
\end{equation}
with
\[
r=-\ell+1+\big|\tvect\big|^2-\big|\svect\big|^2,
\]
where we use the definition from~\cite{DR} for notation \eqref{eq:QNM},
for $r$ possibly negative.
\end{lemma}
\begin{proof}
Let $\Afr$ denote the amalgamated free product von Neumann
algebra on the left hand side of~\eqref{eq:QNM}.
If $N$ is commutative, then the isomorphism~\eqref{eq:QNM} follows directly from~\cite[Cor.\ 3.2]{Dy}.
For general $N$, let $f_j$ be a minimal projection in the $j$th summand of $N$ and let
$\displaystyle f=\sum_{j=1}^\ell f_j$.
Then, as is well known,
\begin{equation}\label{eq:fAf}
f\Afr f\cong(Q\otimes fNf)*_{fNf}fMf.
\end{equation}
(This can be seen directly, or by applying~\cite[Lem.\ 2.2]{BD} and taking the obvious
representation of the full free product C$^*$--algebra.)
The inclusion matrix for \linebreak $fNf\subset fMf$ is still $\Lambda$, while the trace
vectors for $\tau(f)^{-1}\tau |_{fMf}$ and $\tau(f)^{-1}\tau |_{fNf}$
are $\tau(f)^{-1}\svect$ and $\tau(f)^{-1}\tvect$, respectively.
Since $fNf$ is commutative, from~\eqref{eq:fAf} we get
\begin{equation}\label{eq:fAf=}
f\Afr f\cong Q_{\frac{\tau(f)}{t_1}}*\cdots*Q_{\frac{\tau(f)}{t_\ell}}*L(\Fb_{r'}),
\end{equation}
with
\[
r'=-\ell+1+\tau(f)^{-2}(\big|\tvect\big|^2-\big|\svect\big|^2).
\]
Using the formula in~\cite[Prop.\ 4]{DR} for rescalings,
from~\eqref{eq:fAf=} and $\Afr\cong(f\Afr f)_{\tau(f)^{-1}}$, we get~\eqref{eq:QNM}, as required.
\end{proof}

\begin{remark}
\label{remark4.4}
Let $M$ be a II$_{1}$ factor and $\alpha\in\autM$. Then for any projection $p$ in $M$ there exists $\gamma\in\innerM$ such that $\gamma\circ\alpha (p)=p$
\end{remark}
In fact, by the uniqueness of the trace on $M$, $p$ and $\alpha (p)$ must have the same trace. Therefore we can find a unitary $u$ in $M$ such that $u\alpha (p)u^{*}=p$.  

\begin{lemma}
\label{lemma4.5}
Consider a II$_{1}$ factor and an automorphism $\alpha$ of $M$. Suppose that $p$ is a projection in $M$ with $\alpha(p)=p$. Then, $\alpha$ is inner if and only if $\alpha| _{pMp}\in\operatorname{Int}(pMp)$.
\end{lemma}

\begin{proof}
One of the implications is obvious since if $\alpha=\Ad u$ fixes $p$ then $p$ commutes with $u$ and $\alpha |_{pMp}=\Ad (pu)$ on $pMp$. For the other one, let $\tM=M\bar{\otimes}B(\ell ^{2}(\mathbb{Z}))$ be the associated II$_{\infty}$ factor and set $\widetilde{\alpha}=\alpha\otimes Id_{B(\ell ^{2}(\mathbb{Z}))}\in \operatorname {Aut}(\tM)$. It will suffice to show that $\widetilde{\alpha}\in \operatorname {Int}(\tM)$. Indeed, $\widetilde{\alpha}(1\otimes e_{11})=1\otimes e_{11}$, so $\widetilde{\alpha}$ inner implies that $\alpha=\widetilde{\alpha}|_{(1\otimes e_{11})\tM (1\otimes e_{11})}$ is inner.

Let $q=p\otimes 1$. Take partial isometries $v_{j}$ in $\tM$ with $j\geq 0$, such that for every $j$ we have $v_{j}^{*}v_{j}=q$ and $\displaystyle q+\sum _{j=1}^{\infty}v_{j}v_{j}^{*}=1$. By hypothesis, there exists $u\in q\tM q$ such that $\widetilde{\alpha}(x)=uxu^{*}$ for all $x\in q\tM q$. Set 
$$
w=\sum _{i=1}^{\infty}\widetilde{\alpha}(v_{i})uv_{i} ^{*},
$$
which is a unitary in $\tM$. In addition, for any $y\in\tM$ we have  
\begin{eqnarray*}
wyw^{*}&=&\displaystyle\sum _{i,j=1}^{\infty}\widetilde{\alpha}(v_{i})uv_{i} ^{*}yv_{j} ^{*}u\widetilde{\alpha}(v_{j})=\sum _{i,j=1}^{\infty}\widetilde{\alpha}(v_{i})\widetilde{\alpha}(v_{i}^{*}yv_{j})\widetilde{\alpha}(v_{j}^{*})=\\
& & \displaystyle\sum _{i,j=1}^{\infty}\widetilde{\alpha}(v_{i}v_{i} ^{*}yv_{j}v_{j} ^{*})=\widetilde{\alpha}(y).                                               
\end{eqnarray*}
so $\widetilde{\alpha}\in\operatorname{Int}(\tM)$.
\end{proof}

As an easy consequence of the previous lemma we obtain the following corollary.

\begin{corollary}
\label{corollary4.6}
Let $M$ be a II$_{1}$ factor and $\alpha$ an automorphism of $M$. Let $p$ be a projection in $M$ fixed by $\alpha$. Then $\alpha$ and $\alpha|_{pMp}$ have the same outer invariant.
\end{corollary}

Now we are ready to show that for any II$_{1}$ factor $Q$, the factor $\freeprodRRo$ is the enveloping algebra of an inclusion $A_{1}\subset B_{1}$, with $A_{1}$ a free product factor. Observe that if the commuting square  
\[\hspace{4cm}\comsquare\]
is extremal, then the von Neumann algebras $\freeproductDC$ and $\freeproductBA$, are II$_{1}$ factors.

\begin{proposition}
\label{proposition4.7} 
   Let \vspace{-.5\baselineskip}
   \[\hspace{4cm}\comsquare\]
   be the extremal commuting square of Lemma \ref{lemma4.2}. Consider the inclusion of amalgamated free products 
\begin{equation}
\label{freein}
\freeproductDC\subset\freeproductBA ,
\end{equation}
where $Q$ is a $II_{1}$ factor. Set $\theta _{1}=\ad _{[\freeproductDC]}\,(U_{0})$. Then, $\theta_{1}$ is a $\luna$-action on $\freeproductDC$ such that $\freeproductBA$ is isomorphic to the crossed product $(\freeproductDC)\rtimes _{\theta _{1}}\luna$. Moreover, the outer invariant of $\theta _{1}$ is $(p, \obstruction)$, and $\freeproductDC$ is isomorphic to $\left (\ast _{1} ^{p}Q_{\frac{1}{p}}\right )\ast\mathfrak{L}\left (\Fb _{r}\right )$, where $r=-p+1+\frac{1}{p}-\frac{1}{p^3}$.
 
\noindent Lastly, iteration of the Jones basic construction for the inclusion of factors in (\ref{freein}) yields 
   \[ \A=\freeproductRR\cong\freeprodRRo, \]
as the enveloping algebra.
\end{proposition}

\begin{proof}
    First note that the span of elements of the form $q_{1}c_{1}q_{2}c_{2}\cdots q_{n-1}c_{n-1}q_{n}c_{n}$ with $q_{i}\in Q$ and $c_{i}\in C$ is a dense set in $\freeproductDC$, as described in \cite{Popa1}. It follows that  $\ad_{[\freeproductBA ]}(U_{0})$ leaves $\freeproductDC$ invariant, since $U_{0}$ commutes with every element of $Q$ and $\ad _{B}(U_{0})$ leaves $C$ invariant. 

Moreover, using the properties satisfied by the $U_{k}$'s and $g$, it is easy to check that 
\begin{equation}
\label{eqn1}
\ad _{C}(U_{0}^{p})=\ad_{C}(g^{*}).
\end{equation}
Also, since $g$ commutes with the elements of $Q$, (\ref{eqn1}) still holds if we replace $C$ by $\freeproductDC$ , i.e.,
\begin{equation*}
\ad _{[\freeproductDC ]}\,(U_{0}^{p})=\ad _{[\freeproductDC ]}\,(g^{*}).
\end{equation*}
Thus, $\theta _{1}=\ad _{[\freeproductDC ]}\,(U_{0})$ is a $\luna$-action on $A_{1}=\freeproductDC$ with outer invariant $(p,\obstruction)$. 

In addition, by the definition of the von Neumann algebras $A,\, B,\, C$, and $D$ in the commuting square and the representation of any element in the amalgamated free 
product we described before, it follows that
\begin{equation*}
\{(\freeproductDC ) \cup \{U_{0}\}\}^{\prime\prime}=\freeproductBA.
\end{equation*}
Set $B_{1}=\freeproductBA$ and denote by $\tau$ the trace on $B_{1}$, as defined by Popa in \cite{Popa1}. We want to show that $B_{1}$ can be realized from  $A_{1}$ through a crossed product construction with respect to the $\gruppo$-action $\theta _{1}$. Observe that any monomial which is an alternating product of elements in $Q\otimes B$ and $A$ can be written using only one occurrence of $U_{0}$ to some power, because of the relations (ii) and (iii) of Lemma \ref{lemma4.1} and the condition $[Q,U_{0}]=0$. Therefore, it suffices  to show that if $m$ is such a monomial and it contains a non-zero power of $U_{0}$, then $\tau (m)=0$. 

Denote by $E_{1}^{D}$ the conditional expectation defined on $Q\otimes D$ and by $E_{2} ^{D}$ the conditional expectation on $C$. Using again the relations (ii) and (iii) of Lemma \ref{lemma4.1} and the fact that $Q$ and $U_{0}$ commutes, we obtain that it suffices to consider elements of the form 
\begin{equation*}
a=U_{0}^{k}q_{0}c_{0}q_{1}c_{1}\hdots q_{i}c_{i}q_{i+1}c_{i+1}\cdots q_{n}c_{n}q_{n+1}
\end{equation*}
where $q_{i}\in (Q\otimes D)\cap ker E_{1} ^{D}$, $c_{i}\in C\cap  ker E_{2} ^{D}$ $\forall i$. Indeed $m$ is a linear combination of elements of this form. Since $\tau(U_{0}^kc_{1})=0$ (see the definition of the trace on $R_{0}$ in Lemma \ref{lemma4.1}), the definition of the trace  on the amalgamated free product implies that $\tau (m)=0$. 
 
Now denote by $\A$ the enveloping algebra obtained by iteration of the Jones basic construction for the inclusion 
\begin{equation*}
\freeproductDC\subset\freeproductBA .
\end{equation*}
By Lemma \ref{lemma4.2} and Theorem 1.2  in \cite{Radu1} we obtain that 
\begin{equation*}
\A\cong\freeproductRR\cong\freeprodRRo.
\end{equation*}
The last step of the proof is to show that $\freeproductDC$ is isomorphic to 
$\left (\ast _{1} ^{p}Q_{\frac{1}{p}}\right )\ast\mathfrak{L}\left (\Fb _{r}\right )$.
Note that $D=alg\{g\}=\mathbb{C}g_{1}\oplus\mathbb{C}g_{2}\hdots\oplus\mathbb{C}g_{p}\cong\mathbb{C}^{p}$, where the $g_{i}$'s are the spectral projections of $g$. Moreover a maximal set of mutually orthogonal minimal central projections for the finite dimensional algebra $C$ is given by the spectral projections of $U_{-1}^{p}$, since $U_{-1}^{p}$ generates the center of $C$. Thus, by \cite{GHJ} $C$ is isomorphic to $M_{p}(\mathbb{C})\oplus M_{p}(\mathbb{C})\hdots\oplus M_{p}(\mathbb{C})$, where the sum has $p$ terms.  
In addition, the inclusion matrix for the inclusion $D\subset C$ is given by
    
\addvspace{-.4\baselineskip}
\begin{equation*}
\Lambda=\begin{pmatrix}
1 & 1 & 1 & \cdots & 1 \\
1 & 1 & 1 & \cdots & 1 \\
1 & 1 & 1 & \cdots & 1 \\
\vdots & \vdots & \vdots &\cdots &\vdots \\
1 & 1 & 1 & \cdots & 1
\end{pmatrix}
\end{equation*} 

\noindent and the trace vector for the trace on $C$ by

   \[\vec{s}= \frac{1}{p^2} \left (\begin{array}{c}
            1 \\
            1\\
            1 \\
            \vdots \\
            1 
        \end{array}\right ). \]

\vspace{.7\baselineskip}
\noindent By Lemma~\ref{lemma4.3},
\[A_{1}=(Q\otimes D)\ast _{D} C\cong\left (\ast _{1} ^{p}Q_{\frac{1}{p}}\right )\ast\mathfrak{L}\left (\Fb _{-p+1+\frac{1}{p}-\frac{1}{p^3}}\right ).\]

\end{proof}

Taking $A_{1}=\freeproductDC$ and $B_{1}=\freeproductBA$, we can reformulate the previous proposition as follows.

\begin{theorem}
\label{theorem4.8}
Let $N$ be any II$_{1}$ factor and $p$ a prime. Let $Q$ be a II$_1$ factor such that $Q_{t}=N\otimes M_{p}(\mathbb{C})$ for $t=\sqrt{\frac{p+1}{p^3}}$. Then the II$_{1}$ factor $\A=\freeprodRRo$ can be realized as the enveloping algebra of an inclusion $A_{1}\subset B_{1}$, where $A_{1}$ is a free product factor and $B_{1}$ is the crossed product of $A_{1}$ by a $\luna$ action $\theta _{1}$. Moreover, by perturbing $\theta_{1}$ by an inner automorphism, if necessary, we obtain a $\luna$-action on $\freeprodN$, which has outer invariant $(p,\obstruction)$. We denote this action  on $\freeprodN$ by $\widetilde{\theta _{1}}$.  

\end{theorem}
\begin{proof}
We showed in the previous theorem that $A_{1}\cong\left (\ast _{1} ^{p}Q_{\frac{1}{p}}\right )\ast\mathfrak{L}(\Fb _{r})$ with $r=-p+1+\frac{1}{p}-\frac{1}{p^3}$. Denote by $\tau$ the trace on $A_{1}$. Let $q$ a projection in $A_{1}$ such that $\tau (q)=\sqrt{\frac{p+1}{p^3}}=t$. By Remark \ref{remark4.4} there exists $\phi _{1}\in\operatorname{Int}(A_{1})$ such that $\phi _{1}\circ\theta _{1}(q)=q$. Note that $\theta$ and $\phi _{1}\circ\theta  _{1}$ have the same outer invariant since they differ by an inner automorphism. In addition, using the rescaling formula in \cite{DR} we obtain that $\displaystyle qA_{1}q\cong\freeprodN$, since $N=(Q_{\frac{1}{p}})_{t}$. Moreover, Corollary \ref{corollary4.6} implies that the automorphism $\widetilde{\theta _{1}}=(\phi _{1}\circ\theta _{1})|_{qA_{1}q}$ has also outer invariant $(p,\obstruction)$.
\end{proof}

The automorphism $\widetilde{\theta _{1}}$ is our first example of a $\luna$-action on the free group factor $\freeprodN$ with outer invariant $(p, \obstruction )$. The second example example of a $\luna$-action on $\freeprodN$ with the same outer invariant as $\widetilde{\theta _{1}}$ is also obtained through a subfactor construction.
Again we let $Q$ be a II$_1$ factor such that $Q_{t}=N\otimes M_{p}(\mathbb{C})$.
The $\luna$-action $\widetilde{\theta _{2}}$ is constructed using the crossed product factor $\M=\crossprod$, where $\gamma$ is a $\gruppo$-action on $\pippo$. $\gamma$ is defined using two $\gruppo$--kernels, $\alpha\in\operatorname{Aut}\left (\ifreeprod\right )$ and $\beta\in\operatorname{Aut}(R)$, with conjugate obstructions to lifting. 

Let $\displaystyle\lambda_{i}:Q\longrightarrow \ast_{1}^{p}Q$ be the canonical map defined on the i-th copy of $Q$ in the free product, and denote by $u$ the unitary which generates $\mathfrak {L}\left (\Fb _{1-\frac{1}{p}}\right )\cong\mathfrak{L}(\mathbb{Z} _{p})$. We define an automorphism $\alpha$ on $\ifreeprod$ by: 
\begin{itemize}
\item[\hspace{-1cm}]
           $\alpha (\lambda _{i}(x))=\lambda _{i+1}(x)$ for every $x\in Q$, with $i\in\{1,\hdots ,p-1\}$ 
\item[\hspace{-1cm}]
           $\alpha (\lambda _{p}(x))=u\lambda _{1}(x)u^{*}$, for all $x\in Q$
\item[\hspace{-1cm}]
          $\alpha (u)=\obstruction u$.
\end{itemize}
Note that $\alpha ^{p}=\Ad u$, so that $\alpha$ has outer invariant $(p,\obstruction)$.

Using the model  $R\cong R_0=\{g, U_0, U_1,\hdots\}^{\prime\prime}$ for the hyperfinite $II_1$ factor, we define the automorphism $\beta=\ad _{R_0}(U_{-1})$. Since $\beta ^{p} =\ad _{R_0}(g)$, $\beta$ is a $\gruppo$--kernel with obstruction $\obstruconj$ to lifting. Observe that $\alpha\otimes\beta$ has outer invariant $(p, 1)$ so it can be perturbed by an inner automorphism to obtain an action of $\gruppo$ on $\tenpro$ (see comment at the end of Section 2). Denote by $e_{i}$, for $i=1,\hdots p$, the spectral projections of the unitary $u$ so that $u=\displaystyle\sum _{j=1} ^{p}e^{\frac{2\pi ij}{p}}e_{j}$. Then we can define the action on $\tenpro$ by
\begin{equation*}
\gamma=\left (\ad _{\left [\tenpro\right ]} W \right)\alpha\otimes\beta ,
\end{equation*}
where $W=\displaystyle\sum _{l=0} ^{p-1}{e^{\frac{-2\pi i l}{p^2}}E_{l}}$ with $\displaystyle E_{l}=\sum _{\substack {k,j=1,\hdots ,p\\ k+j\equiv l\bmod p}}{e_{k}\otimes g_{j}}$, is a p-root of $u^{*}\otimes g^{*}$ which  belongs to the fixed point algebra of $\alpha\otimes\beta$ in $\tenpro$.

\begin{theorem}
\label{theorem4.9}
Let $\M=\crossprod$ and denote by $v$ the unitary implementing the crossed product.  
Set 
\[A_{2}=((\ast _{i}^{p}Q\otimes 1)\cup \{u\otimes 1, 1\otimes g, v\})^{\prime\prime},\] 
where $[u\otimes 1,1\otimes g]=0$, $v(u\otimes 1)=\obstruction (u\otimes 1)v$, and $v(1\otimes g)=\obstruconj (1\otimes g)v$.

Consider the inclusion of subalgebras of $\M$
\[A_{2}\subset B_{2}=(A_{2}\cup \{1\otimes U_{0}\})^{\prime\prime}.\]
Then, $A_{2}$ is isomorphic to  $\left (\ast _{1} ^{p}Q_{\frac{1}{p}}\right )\ast\mathfrak{L}(\Fb _{r})$ with $r=-p+1+\frac{1}{p}-\frac{1}{p^3}$. Also, $ B_{2}$ is equal to the crossed product $A_{2}\rtimes _{\theta _{2}}\luna$, where $\theta _{2}=\Ad (1\otimes U_{0})$ is a $\luna$-action on $A_{2}$ with outer invariant $(p, \obstruction)$. 
  
\noindent Furthermore, $\M$ is the enveloping algebra for the inclusion $A_{2}\subset B_{2}$. 
\end{theorem}

\begin{proof} 
 First observe that if $\{g_{k}\mid 1\leq k\leq p\}$ denotes the set of spectral projections of $g$ as in Lemma \ref{lemma4.1}, then $\Ad (1\otimes U_{0})(1\otimes g_{k})=1\otimes g_{k+1}$ for $1\leq k< p$ and $\Ad (1\otimes U_{0})(1\otimes g_{p})=1\otimes g_{1}$. Also, $\Ad (1\otimes U_{0})$ leaves $u\otimes 1$ invariant. Therefore, if $E_{l}=\displaystyle\sum _{\substack{k,j=1,\hdots ,p\\k+j\equiv l\bmod p}}{e_{k}\otimes g_{j}}$ and $W=\displaystyle\sum _{l=0}^{p-1}{e^{\frac{-2\pi i l}{p^2}}E_{l}}$ as before, we have that $\Ad (1\otimes U_{0})(E_{l})=E_{l+1}$ for $0\leq l\leq p-1$ and $\Ad (1\otimes U_{0})(E_{p-1})=E_{0}$. Consequently,
\begin{eqnarray*}
& \hspace{-.7cm}& (1\otimes U_{0})v(1\otimes U_{0})^{*}= (1\otimes U_{0})\gamma (1\otimes U_{0}^{*})v=(1\otimes U_{0})\Ad W(\beta (1\otimes U_{0}^{*})v=\\
&\hspace{-.7cm} &e^{\frac{-2\pi i}{p^2}}(1\otimes U_{0})\Ad W(1\otimes U_{0}^{*})v=e^{\frac{-2\pi i}{p^2}}\Ad (1\otimes U_{0})(W)W^{*}v=\\
&\hspace{-.7cm} &e^{\frac{-2\pi i}{p^2}}\left (\sum _{l=0}^ {p-2}e^{\frac{-2\pi il}{p^2}}E_{l+1}+e^{\frac{2\pi i(1-p)}{p^2}}E_{0}\right )W^{*}v=\left (\sum _{l=1}^{p-1}E_{l}+\obstruconj E_{0}\right )v
\end{eqnarray*}

Moreover, if $\displaystyle \lambda _{i}:Q\longrightarrow \ast_{i}^{p}Q$ is the canonical map on the i-th copy of $Q$ in the free product, then $\Ad (1\otimes U_{0})(\lambda _{i}(x)\otimes 1)=\lambda_{i}(x)\otimes 1$ for all $x\in Q$ and for every $1\leq i\leq p$, thus $\theta _{2}=\Ad (1\otimes U_{0})$ defines an outer automorphism of $A_{2}$. Also, $\ad _{A_{2}}(1\otimes U_{0})^{p}=\ad _{A_{2}}(1\otimes g^{*})$ and $\ad _{A_{2}}(1\otimes U_{0})(1\otimes g^{*})=\obstruction (1\otimes g^{*})$, so $\theta _{2}$ has outer invariant $(p, \obstruction )$.

Note that any monomial $m$ in $B_{2}$ can be written as $m=m^{\prime}v^{k}$, with $m'$ a monomial in $\displaystyle ((\ast_{1}^{p}Q\otimes 1)\cup\{u\otimes 1, 1\otimes g, 1\otimes U_{0}\})^{\prime\prime}$ and $0\leq k\leq p-1$. Clearly $m'$ can be written using only one occurrence of $1\otimes U_{0}$ to some power so the same holds for $m$.

To complete the proof that $B=(A\cup\{1\otimes U_0\})^{\prime\prime}$ is the crossed product of $A_{2}$ by $\theta _{2}=\ad _{A_{2}}(1\otimes U_0)$, we need to show that the trace $\tau$, which $B_{2}$ inherited as a subalgebra of $\M=\crossprod$, is zero on any monomial $m$ in $B_{2}$ containing $1\otimes U_{0}$. If the monomial has the form $m=m^{\prime}v^{k}$ with  $m^{\prime}$ a monomial in $\displaystyle (*_{1}^{p}Q\otimes 1)\cup\{u\otimes 1, 1\otimes g, 1\otimes U_{0}\mid 1\leq i\leq p\}$ containing $1\otimes U_{0}$, and $0< k\leq p-1$. Since $v$ is the unitary implementing the crossed product, we can conclude immediately that $m$ has zero trace. Next, assume that $k=0$, so that the trace of $m$ as an element of $\M$ coincides with its trace as an element of $\tenpro$. By the definition of the trace on $R_{0}$ (see Lemma \ref{lemma4.1}) it follows that $\tau (m)=0$.
    
Observe also that as in Lemma \ref{lemma4.1}, the successive steps of the Jones basic construction for the inclusion $A_{2}\subset B_{2}$ are obtained by adding the unitaries $U_{1}, U_{2},\hdots$. Therefore, $\M$ is the enveloping algebra for $A_{2}\subset B_{2}$.     

Lastly, we want to show that $A_{2}$ is an interpolated free group factor. Using the spectral projections $\{1\otimes g_{i}\}_{1} ^{p}$ of $1\otimes g$, and the unitary $v$ we can define the family of partial isometries
$$
f_{j,1} =(1\otimes g_{j})v^{j-1}(1\otimes g_{1})=(1\otimes g_{j})v^{j-1}
$$
Consider the induced algebra $f_{1,1}A_{2}f_{1,1}$. By Lemma 5.21 in \cite{VDN},  \[ f_{1,1}A_{2}f_{1,1}=\{f_{j,1}^{*}xf_{k,1}\mid 1\leq k,j\leq p, x\in\Omega\}^{\prime\prime },\]
 where $\Omega =\{\lambda _{h}(q)\otimes 1, u\otimes 1, 1\otimes g, v\mid 1\leq h\leq p \text{ and }$q$ \text{ is a generator of }Q\}$ is the set of generators for $A_{2}$. We claim that 
\[f_{1,1}A_{2}f_{1,1}=\left (\ifreeprod\right )\otimes g_{1}.\] 
This is a simple consequence of computing the various products of the form $f_{j,1}^{*}xf_{k,1}$ for $x\in\Omega$. Indeed these calculations yields
\begin{itemize}
\item[\hspace{-1cm}i)] $f_{j,1}^{*}(1\otimes g)f_{j,1}=e^{\frac{2\pi i j}{p}}(1\otimes g_{1})$, while $f_{j,1}^{*}(1\otimes g)f_{k,1}=0$ for $j\neq k$ and $1\leq j,k\leq p$  
\item[\hspace{-1cm}ii)] $f_{j,1}^{*}(u\otimes 1)f_{j,1}=e^{\frac{-2\pi i (j-1)}{p}}(u\otimes g_{1})$, while $f_{j,1}^{*}(u\otimes 1)f_{k,1}=0$ for $j\neq k$ and $1\leq j,k\leq p$ 
\item[\hspace{-1cm}iii)] $f_{j+1,1}^{*}vf_{j,1}=1\otimes g_{1}$ for $1\leq j\leq p$ 
\item[\hspace{-1cm}iv)] $f_{j,1}^{*}vf_{k,1}=0$ for $j\neq k+1$ and $1\leq j,k\leq p$, except for $f_{1,1}^{*}vf_{p,1}=1\otimes g_{1}$. 
\end{itemize}
Also, if $\{e_{j}\}_{1} ^{p}$ denote the spectral projections of $u$, then for any $1\leq h\leq p$ we have that
\begin{itemize}
\item[\hspace{-1cm}v)] $\displaystyle f_{j,1}^{*}(\lambda _{h}(q)\otimes 1)f_{j,1}=\sum _{l,k=1} ^{p} e^{\frac{2\pi i (l-k)(j-1)}{p^2}}e_{l}\lambda _{h-j+1}(q)e_{k}\otimes g_{1}$, for all $\lambda _{h}(q)\otimes 1\in \Omega$ and $1\leq j< h+1$.  
\item[\hspace{-1cm}vi)]  $\displaystyle f_{j,1}^{*}(\lambda _{h}(q)\otimes 1)f_{j,1}=\sum _{l,k=1} ^{p} e^{\frac{2\pi i (l-k)(j-1-p)}{p^2}}e_{l}\lambda _{p-j+1+h}(q)e_{k}\otimes g_{1}$, for all $\lambda _{h}(q)\otimes 1\in \Omega$ and $h+1\leq j\leq p$.  
\item[\hspace{-1cm}vii)] $f_{j,1}^{*}(\lambda _{h}(q)\otimes 1)f_{k,1}=0$ for $j\neq k$ and $1\leq j,k\leq p$. 
\end{itemize}

Therefore, $f_{1,1}A_{2}f_{1,1}\cong \ifreeprod$. Using the scaling formula for interpolated free group factors we conclude that $A_{2}\cong \left (\ast _{1} ^{p}Q_{\frac{1}{p}}\right )\ast\mathfrak{L}(\Fb _{r})$ with $r=-p+1+\frac{1}{p}-\frac{1}{p^3}$. 
\end{proof} 

By compressing the II$_{1}$ factor $A_{2}$ and using the rescaling formula in \cite{DR} we obtain a $\luna$-action on a free product factor, with the same outer invariant as $widetilde{\theta_{1}}$.

 \begin{corollary}
\label{corollary4.10}
Let $N$ be any II$_{1}$ factor. Take the II$_{1}$ $Q$ so that $Q_{t}=N\otimes M_{p}(\mathbb{C})$ for $t=\sqrt{\frac{p+1}{p^3}}$. Then the II$_{1}$ factor $\M=\crossprod$ is the enveloping algebra of an inclusion $A_{2}\subset B_{2}$, where $A_{2}$ is a free product factor, and $B_{2}$ is the crossed product of $A_{2}$ by a $\luna$ action $\theta _{2}$. Let $q$ be a projection of trace $t=\sqrt{\frac{p+1}{p^3}}$ and $\phi _{2}$ an inner automorphism such that $\phi _{2}\circ\theta_{2}(q)=q$ (as in Remark \ref{remark4.4}). Then $\widetilde{\theta} _{2}=(\phi _{2}\circ\theta _{2})|_{qA_{2}q}$ is a $\luna$-action on $\freeprodN$, which has outer invariant $(p,\obstruction)$.
\end{corollary}

\section{The Connes invariant of the factors $\A$ and $\M$}  
We devote this section to the computation of the Connes invariant of the factor $\A=\freeprodRRo$, where $R_{0}=R_{-1}\rtimes _{\Ad (U_{0})}\luna$. We show in Corollary \ref{corollary5.2} that this is equivalent to computing the Connes relative invariant of the inclusion $R_{-1}\subset R_{0}$. The proof is based on an adaption of the method of Sakai's Lemma 4.3.3 in \cite{Sakai}, to amalgamated free products of von Neumann algebras.  

Suppose $\Nc$ is a von Neumann algebra that is unitally embedded in von Neumann algebras
$A_1$ and $A_2$, each equipped with a faithful, normal tracial state $\tau_1$, respectively $\tau_2$,
such that the restrictions of $\tau_1$ and $\tau_2$ to $\Nc$ agree.
Let $E_i:A_i\to\Nc$ be the $\tau_i$--preserving conditional expectation.
Let
\[
(\Mcal,E)=(A_1,E_1)*_\Nc(A_2,E_2)
\]
be the amalgamated free product of von Neumann algebras.
Then, by \cite[Proposition 3.1]{Popa2} $\tau_1\circ E$ is a normal, faithful tracial state on $\Mcal$.
As is conventional, we will denote the map $\Mcal\to L^2(\Mcal,\tau)$ by $x\mapsto\xh$.
For emphasis, we may denote the left action of $\Mcal$ on $L^2(\Mcal,\tau)$ by $x\mapsto\lambda(x)$,
and we will denote the right action of $\Mcal^{\operatorname{op}}$
on $L^2(\Mcal,\tau)$ by $x\mapsto\rho(x)=J\lambda(x)J$,
where $J\xh=\widehat{x^*}$.

\begin{proposition}
\label{proposition5.1}
Suppose there are unitaries $u_1\in A_1$ and $u_2,u_3\in A_2$ such that
\begin{gather}
x\in A_1,\;E_1(x)=0\implies E_1(u_1xu_1^*)=0 \label{eq:Eu1} \\
E_1(u_1)=0=E_2(u_2)=E_2(u_3)=E_2(u_2^*u_3). \label{eq:uorthog}
\end{gather}
Let $y\in\Mcal$ and set 
\begin{equation}\label{eq:eps}
\eps=\max_{1\le i\le3}\|[y,u_i]\|_2,
\end{equation}
where $[a,b]=ab-ba$.
Then
\[
\|y-E(y)\|_2\le14\eps.
\]
\end{proposition}
\begin{proof}
Let $A_i\oup=\ker E_i\subseteq A_i$.
From the construction of the free product with amalgamation, due to Voiculescu~\cite{V85},
$\Mcal$ has subspaces $\Nc$ and
\begin{equation}\label{eq:Asubs}
(A_{i_1}\oup A_{i_2}\oup\cdots A_{i_n}\oup)_{n\ge1,\,i_1,\ldots,i_n\in\{1,2\},\,i_j\ne i_{j+1}},
\end{equation}
which are mutually orthogonal with respect to $E$, hence also with respect to $\tau$,
and $\Mcal$ is the weak closure of the linear span of the union of these subspaces.
Let
\[
F=A_1\oup+A_2\oup A_1\oup+A_1\oup A_2\oup A_1\oup+A_2\oup A_1\oup A_2\oup A_1\oup+\cdots\subseteq\Mcal
\]
be the span of the union of the these subspaces that end in $A_1\oup$,
and let $\Fc$ denote the closure of this space in $L^2(\Mcal,\tau)$.

Let us first observe that
\begin{equation}\label{eq:u_1F}
\lambda(u_1)\rho(u_1)\Fc+\Fc=L^2(\Mcal,\tau)\ominus L^2(\Nc,\tau).
\end{equation}
Indeed, since $E(u_1)=0$, it is easy to see that all of the subspaces~\eqref{eq:Asubs}
belong to either $F$ or $u_1Fu_1^*$, which become $\Fc$ and $\lambda(u_1)\rho(u_1)\Fc$
upon embedding in $L^2(\Mcal,\tau)$.
For example, $A_1\oup A_2\oup=u_1(u_1^*A_1\oup A_2\oup u_1)u_1^*$
and
\[
u_1^*A_1\oup A_2\oup u_1\subseteq A_1\oup A_2\oup u_1+A_2\oup u_1
\subseteq A_1\oup A_2\oup A_1\oup+A_2\oup A_1\oup\subseteq F.
\]
Furthermore, from~\eqref{eq:Eu1}, we find $u_1A_1\oup u_1^*=A_1\oup$, and from 
this we see $\lambda(u_1)\rho(u_1)\Fc\subseteq L^2(\Mcal,\tau)\ominus L^2(\Nc,\tau)$.

Let us now show that the three subspaces
\begin{equation}\label{eq:u_23F}
\Fc,\quad\lambda(u_2)\rho(u_2)\Fc,\quad\lambda(u_3)\rho(u_3)\Fc
\end{equation}
are contained in $L^2(\Mcal,\tau)\ominus L^2(\Nc,\tau)$ and are pairwise orthogonal.
We already have $\Fc\subseteq L^2(\Mcal,\tau)\ominus L^2(\Nc,\tau)$.
Since $E(u_2)=E(u_3)=0$, it is easy to see that for $i=2,3$, we have
\begin{align*}
u_iFu_i^*&\subseteq A_1\oup u_i^*+A_2\oup A_1\oup u_i^*+A_1\oup A_2\oup A_1\oup u_i^*+\cdots \\
&\subseteq A_1\oup A_2\oup+A_2\oup A_1\oup A_2\oup+A_1\oup A_2\oup A_1\oup A_2\oup+\cdots
\end{align*}
Therefore, $\lambda(u_i)\rho(u_i)\Fc\subseteq L^2(\Mcal,\tau)\ominus L^2(\Nc,\tau)$ for $i=2,3$ and
both are orthogonal to $\Fc$.
It remains to show that these two subspaces are orthogonal to each other.
However, letting
\[
a_j,a_j'\in\begin{cases}
A_1\oup&j\text{ odd,} \\
A_2\oup&j\text{ even}
\end{cases}
\]
using freeness we have
\begin{gather*}
\langle a_na_{n-1}\cdots a_2a_1u_2^*,a_n'a_{n-1}'\cdots a_2'a_1'u_3^*\rangle_\tau
=\tau(u_2a_1^*a_2^*\cdots a_n^*a_n'\cdots a_2'a_1'u_3^*) \\
=\tau(u_2E(a_1^*E(a_2^*\cdots E(a_{n-1}^*E(a_n^*a_n')a_{n-1}')\cdots a_2')a_1')u_3^*) \\
=\tau(E(\cdots)u_3^*u_2)
=\tau\circ E(E(\cdots)u_3^*u_2)=\tau(E(\cdots)E(u_3^*u_2))=0,
\end{gather*}
where we use $E(u_3^*u_2)=0$.
A calculation shows
\[
\langle a_na_{n-1}\cdots a_2a_1u_2^*,a_m'a_{m-1}'\cdots a_2'a_1'u_3^*\rangle_\tau=0\text{ when }n\ne m. 
\]
These calculations imply the orthogonality of $\lambda(u_2)\rho(u_2)\Fc$ and $\lambda(u_3)\rho(u_3)\Fc$.

Let $P_0$ be the projection from $L^2(\Mcal,\tau)$ onto $\Fc$ and for $i=1,2,3$ let
$P_i$ be the projection from $L^2(\Mcal,\tau)$ onto $\lambda(u_i)\rho(u_i)\Fc$.
Since $\lambda(u_i)$ and $\rho(u_i)$ are unitaries, we have
\begin{equation}\label{eq:Pi}
P_i=\lambda(u_i)\rho(u_i)P_0\rho(u_i^*)\lambda(u_i^*),\qquad(1\le i\le3).
\end{equation}
Note that $\|\lambda(u_i)\zeta\|_2=\|\rho(u_i)\zeta\|_2=\|\zeta\|_2$ for all $\zeta\in L^2(\Mcal,\tau)$,
since this holds when $\zeta=\xh$ for $x\in\Mcal$.

Let $y\in\Mcal$ and let $\eps$ be as in~\eqref{eq:eps}.
For $1\le i\le 3$, using~\eqref{eq:Pi} we get
\begin{align*}
\big|\|P_i(\yh)\|_2-\|P_0(\yh)\|_2\big|
&=\big|\|P_0(\widehat{u_iyu_i^*})\|_2-\|P_0(\yh)\|_2\big|
\le\|P_0(\widehat{u_iyu_i^*}-\yh)\|_2 \\
&\le\|u_iyu_i^*-y\|_2=\|[u_i,y]\|_2\le\eps.
\end{align*}

Let $s=\|y-E(y)\|_2$.
For $1\le i\le3$, we have that $\|P_i(\yh)\|_2\leq s$ and 
\[
\big|\|P_i(\yh)\|_2^2-\|P_0(\yh)\|_2^2\big|
=\big|\|P_i(\yh)\|_2-\|P_0(\yh)\|_2\big|
\big(\|P_i(\yh)\|_2+\|P_0(\yh)\|_2\big)
\le2s\eps.
\]
From~\eqref{eq:u_1F}, we have
\[
s^2\le\|P_0(\yh)\|_2^2+\|P_1(\yh)\|_2^2\le2\|P_0(\yh)\|_2^2+2s\eps,
\]
so $\|P_0(\yh)\|_2^2\ge\frac{s^2}2-s\eps$.
On the other hand, from the orthogonality of the three subspaces~\eqref{eq:u_23F} and their
containment in $L^2(\Mcal,\tau)\ominus L^2(\Nc,\tau)$,
we get
\[
s^2\ge\|P_0(\yh)\|_2^2+\|P_2(\yh)\|_2^2+\|P_3(\yh)\|_2^2
\ge3\|P_0(\yh)\|_2^2-4s\eps\ge\frac32s^2-7s\eps.
\]
This yields immediately $14\eps\ge s$.
\end{proof}

\begin{corollary}
\label{corollary5.2}
Consider a II$_{1}$ factor $N$ and let $K$ a discrete group acting on $N$. For any II$_{1}$ factor $Q$, let $M=(Q\otimes N)* _{N}(N\rtimes _{\psi}K)$. Then the Connes invariant of $M$ is isomorphic to the relative Connes invariant of the pair $N\subset N\rtimes _{\psi}K$.
\end{corollary}

\begin{proof}
 Consider the map $\Phi :\mbox{Aut}(N\rtimes _{\psi}K,N)\rightarrow\mbox {Aut}(M)$, defined by 
    \[\Phi (\alpha)=(Id_{Q}\otimes\alpha |_{N})*_{N}\alpha,\]
A consequence of the previous proposition is that the argument used by R\u{a}dulescu in \cite[Proposition 3]{Radu2} can also be applied to our factor $M$ so $\phi$ is a continuous map sending $\mbox{Int}(N\rtimes _{\psi}K,N)$ into $\mbox {Int}(M)$, and  $\mbox{Ct}(N\rtimes _{\psi}K,N)$ into $\mbox {Ct}(M)$. Note also that in Proposition 3 the hypothesis that the inclusion $N\subset N\rtimes _{\psi}K$ is hyperfinite is not necessary.
\end{proof}

Thus, for the $II_{1}$ factor $\A=\freeprodRRo$, where $Q$ is a II$_{1}$ factor, we have that $\chi (\A)\cong\chi (R_{0},R_{-1})$. To compute the relative Connes invariant we first need to determine  how the algebra $R_{-1} ^{\omega}\cap R_{0} ^{\prime}$ sits inside the $\omega$-central algebra $(R_{0})_{\omega}=R_{0} ^{\omega}\cap R_{0}^{\prime}$.
 
Given a von Neumann algebra $M$ and a subgroup $G$ of $\autM$, denote by $M^{G}$ the algebra containing all the elements of $M$ fixed by $G$.  

\begin{lemma}
\label{lemma5.3}
   Let $R_{-1}=\{gU_{0}^{p}, U_{1}, U_{2},\hdots\}^{\prime\prime}\subset\{g, 
   U_{0}, U_{1}, U_{2},\hdots\}^{\prime\prime}= R_{0}$ be the inclusion of 
   type $II _{1}$ factor described in Lemma \ref{lemma4.1}. Set $\beta=
   \ad _{R_{0}}(U_{-1})$ and denote by $\beta _{\omega}=(\ad _{R_{0}}(U_{-1})) _{\omega}$. Then, $\beta _{\omega}$ is a non-trivial outer action of order $p$, with the property that 
   \[(R_{0}^{\omega}\cap R_{0}^{\prime}) ^{\langle\beta _{\omega }\rangle}= R_{-1}^{\omega}\cap R_{0}^{\prime}. \]
Moreover, $R_{-1}^{\omega}\cap R_{0}^{\prime}\subset R_{0}^{\omega}\cap R_{0}^{\prime}$ is an inclusion of $II_{1}$ factors of finite index.
\end{lemma}
\begin{proof}

Using the relations satisfied by the $U_{k}$'s and $g$, it is easy to check that\linebreak  $\beta _{\omega}^{p}=\ad _{R_{0}}(U_{-1} ^{p})=\ad_{R_{0}}(g)$. Since $g$ belongs to $R_{0}$, $\ad _{R_{0}}(g)$ is a centrally trivial automorphism, which implies that $\beta _{\omega}^{p}=1$ by Remark \ref{remark2.5}. Also, $\beta _{\omega}=(\ad _{R_{0}}(U_{-1}))_{\omega}$ is non-trivial. Indeed,  if $\beta _{\omega}$ was trivial then $\beta$ would belong to $\operatorname{Ct}(R_{0})$ by remark mentioned above. But for the hyperfinite $II_{1}$ factor $R_{0}$ we have that $\operatorname{Ct}(R_{0})=\operatorname{Int}(R_{0})$ by \cite[Corollary 4]{Connes1}, and $\beta=\ad _{R_{0}}(U_{-1})\not\in\operatorname{Int}(R_{0})$. Therefore $\beta _{\omega}$ is a non-trivial $\gruppo$-action.

To show that $(R_{0}^{\omega}\cap R_{0}^{\prime})^{\langle\beta _{\omega }\rangle}= R_{-1}^{\omega}\cap R_{0}^{\prime}$, first note that $R_{-1}$ is the fixed point algebra $R_{0}^{G}$, where $G=\langle\beta\rangle$. In fact, write any element $x$ of $R_{0}=R_{-1}\rtimes _{\ad _{R_{-1}}(U_{0})}\gruppo$ as $\displaystyle\sum _{i=0}^{p-1}a_{i}U_{0}^{i}$ with $a_{i}\in R_{-1}$. Now observe that $\beta (U_{0})=\jonesinv U_{0}$, while $\beta (x)=x$ for any $x\in R_{-1}$. Thus, $\beta _{\omega}$ acts as the identity on $R_{-1}^{\omega}\cap R_{0}^{\prime}$ and $R_{-1}^{\omega}\cap R_{0}^{\prime}\subset (R_{0}^{\omega}\cap R_{0}^{\prime})^{\langle\beta _{\omega }\rangle}$.

On the other hand, from the observation that
\[ (R_{0}^{\omega}\cap R_{0}^{\prime})^{\langle \beta  _{\omega}\rangle}=\left \{x+\beta (x)+\beta ^{2}(x)+\hdots +\beta ^{p^2-1}(x)
   \, |\, x\in R_{0}\right \}^{\omega}\cap R_{0}^{\prime}, \] 
and $x+\beta (x)+\hdots +\beta ^{p^2-1}(x)\in R_{0}^{\langle\beta\rangle}=R_{-1}$, it follows that $\left (R_{0}^{\omega}\cap R_{0}^{\prime}\right )^{\langle\beta _{\omega}\rangle}\subset R_{-1}^{\omega}\cap R_{0}^{\prime}$.

In addition, $\beta _{\omega}$ is outer since if $\beta _{\omega}=\ad _{R_{0}^{\omega}\cap R_{0}^{\prime}}(t)$ for some $t\in R_{0}^{\omega}\cap R_{0}^{\prime}$, then $txt^{*}=\beta _{\omega}(x)=x$\, $\forall\;x\in R_{-1}^{\omega}\cap R_{0}^{\prime}$. This implies that $t\in (R_{-1}^{\omega}\cap R_{0}^{\prime})^{\prime}\cap (R_{0}^{\omega}\cap R_{0}^{\prime})$. But by \cite{Ocneanu2} (see Theorem 15.30 in \cite{EvK} for a proof and note that the hypothesis of a trivial relative commutant for the original subfactor is unnecessary) the relative commutant of $R_{-1}^{\omega}\cap R_{0}^{\prime}$ in $R_{0}^{\omega}\cap R_{0}^{\prime}$ is trivial. Observe also that  $\beta _{\omega}$ has the smallest period in its outer conjugacy class. Therefore $R_{-1}^{\omega}\cap R_{0}^{\prime}$ is a factor by \cite[Theorem 2.5]{Connes2}. In addition, $ R_{0}^{\omega}\cap R_{0}^{\prime}$ is a $II_{1}$ factor by Theorem 15.15 in \cite{EvK}. Also, the index of the inclusion $R_{-1}^{\omega}\cap R_{0}^{\prime} \subset R_{0}^{\omega}\cap R_{0}^{\prime}$ is finite because of the remark below. 
\end{proof}

The fact that $\beta _{\omega}$ is outer provides us with some additional information about the inclusion  $R_{-1}^{\omega}\cap R_{0}^{\prime} \subset R_{0}^{\omega}\cap R_{0}^{\prime}$.

\begin{remark}
\label{remark5.4}
   Let $\alpha_{0}$ an outer action on a finite von Neumann algebra $M$ and set $G=\langle\alpha _{0}\rangle$. In [Example 3, \cite{Wasse}] A. Wasserman showed that the Jones tower for the inclusion $M^{G}\subset M$ is given by 
\[ M^{G}\subset M\subset M\rtimes G\subset M\otimes B(L^{2}(G))\subset\cdots \cdots.\]
This implies that the inclusion \mbox{$M^{G}\subset M$} has the same form as the inclusion \linebreak $M\subset M\rtimes G$, i.e., $M\cong  M^{G}\rtimes G$. It follows that if $\beta _{\omega}$ is the outer action defined in the previous lemma then $R_{0}^{\omega}\cap R_{0}^{\prime}\cong (R_{-1}^{\omega}\cap R_{0}^{\prime})\rtimes _{\beta _{\omega}}\gruppo$. Therefore, the only non-trivial automorphisms of $R_{0}^{\omega}\cap R_{0}^{\prime}$ which are the identity on $R_{-1}^{\omega}\cap R_{0}^{\prime}$ are $\{\beta _{\omega}^{k}\mid 1\leq k\leq p-1\}$.  
\end{remark} 

Given an  inclusion $N\subseteq M$, the normalizer of $N$ in $M$ is given by $\mathcal{N}_{M}(N)=\{u\in R_{0}|\, \ad u(N)=N\}$. The description given in the previous remark for the inclusion $R_{-1}^{\omega}\cap R_{0}^{\prime}\subset R_{0}^{\omega}\cap R_{0}^{\prime}$ allows us to give a precise characterization of the automorphisms in $\CtRoR$.

\begin{corollary}
\label{corollary5.5}  
Let $\alpha$ be an automorphism in $\CtRoR$. Then, there exists $j\in\{1,\hdots ,p\}$ and $u$ in $\normalizer$ such that  
\[ \alpha =\ad _{R_{0}}(U_{-1}^{j}u).\]
Moreover, up to automorphisms in $\innerRoR$, a set of generators for $\normalizer$ is given by 
\begin{equation*}
   \{\ad _{R_{0}}(U_{0}^{j})\mid j=1,\hdots ,8\}\qquad\text{and }
\end{equation*}
\begin{equation*}
   \big\{\ad _{R_{0}}\big (\sum _{j=1} ^{p}\alpha  _{j}g_{j}\big)\big |\alpha _{j}\in\mathbb{T}\text{ for }1\leq j\leq p\big \},
\end{equation*} 
 where the $g_{i}$'s are the spectral projections of $g\in R_{0}$.
\end{corollary}
\begin{proof}
The statement follows from the previous remark and the same sort of arguments used by R\u{a}dulescu in \cite[Corollary 11 and Lemma 12]{Radu2} for the case $p=2$.
\end{proof}

Recall that the Connes relative invariant was defined as 
\[\chi (R_{0}, R_{-1})=\frac{\CtRoR\cap\chiusoRoR}{\innerRoR}.\]
Thus, now that we know the form of any centrally trivial automorphism,  we need to recognize which of these automorphisms are also approximately inner in order to compute $\chi (R_{0}, R_{-1})$. For this purpose we define the Loi invariant of a II$_{1}$ subfactor.  

Let $N\subset M$ be an inclusion of type $II_{1}$ factors with finite index and Jones' tower $N\subset M\subset M_{1}\subset M_{2}\subset
\cdots$. Consider the group ${\mathcal G}$ formed by all families  $\{\alpha _{k}\}_{k\geq 0}$ of trace preserving automorphisms 
$\alpha _{k}$ on $N^{\prime}\cap M_{k}$, which preserve the inclusion $M^{\prime}\cap M_{k}\subset N^{\prime}\cap M_{k}$, and satisfy the 
conditions
\begin{itemize} 
\item[i)] $\alpha _{k+1}$ is an extension of $\alpha _{k}$ for ever $k\geq 0$.
\item[ii)] $\alpha _{k}$ satisfies $\alpha _{k}(e_{j})=e_{j}$ for all Jones projections $e_{j}$ and $j=1,\hdots ,k$. 
\end{itemize}

There is a natural way to produce an element of ${\mathcal G}$ from an automorphism $\alpha$ in $\autMN$. Let $N\subset M$ be as above and  consider an automorphism $\alpha$ of $M$ leaving $N$ invariant. Since $\alpha$ commutes with the conditional expectation $E_{N}$ onto $N$, it is possible to extend $\alpha$ to an automorphism of the algebra $M_{k}$ defining $\alpha (e_{k})=e_{k}$. We still denote by $\alpha$ the extension of the automorphism to $M_{k}$. Consider now the restriction $\alpha _{k}$ of $\alpha$ to the higher relative commutant $N^{\prime}\cap M_{k}$ (which is an automorphism since $\alpha (N^{\prime}\cap M_{k})=N^{\prime}\cap M_{k}$). The family $\{\alpha _{k}\}_{k\geq 0}$ is an element of ${\mathcal G}$ and is known as the Loi invariant of the automorphism $\alpha$. 

\begin{remark}
\label{remark5.6}
Let $N\subset M$ be a strongly amenable, extremal inclusion of factors in the hyperfinite $II_{1}$ factor. Then \cite[Theorem 15.5]{EvK})
\[\alpha\in\chiusoMN\quad\text{if and only if}\quad\alpha _{k}\mid_{N^{\prime}\cap M_{k}}=Id\,\text{ for all }k\in\mathbb{N}.\]
\end{remark}

The Loi invariant provides us with a tool to recognize which centrally trivial automorphisms are also approximately inner, and thus compute the relative Connes invariant of the inclusion $R_{-1}\subset R_{0}$. The argument used in the following proof is similar to the one used by R\u{a}dulescu in \cite{Radu2}. We present it here for the sake of completeness.

\begin{proposition}
\label{proposition5.7}
The group $\chi (R_{0}, R_{-1})$ is generated, modulo $\innerRoR$, by 
\[ Ad_{R_{0}}(U_{-1}^{p-1}U_{0}). \]

\end{proposition}
\begin{proof}
Let $\theta=\varepsilon(\alpha )$ be any element in $\chi (R_{0},R_{-1})$. Then, $\alpha _{k}\mid_{R_{-1}^{\prime}\cap R_{k}}=Id$ for all $k\geq 0$. By Corollary \ref{corollary5.5} any automorphism in $\CtRoR$ is generated by elements of the form 
$\ad _{R_{0}}(U_{-1} ^{k}u)$, with $u=U_{0}^{j}$ for some $j\geq 0$, or $\displaystyle u=\sum _{i=0} ^{p}\alpha_{i}g_{i}$, where $\alpha _{i}\in\mathbb{T}$ for every $j\geq 0$ and $\alpha _{i}\in\mathbb{T}$. We want to show that $\ad _{R_{0}}(U_{-1} ^{k}u)$ belongs to $\chiusoRoR$ if and only if $u=U_{0}^{j}$ and $k+j\equiv 0\bmod p$.
 
Consider the automorphism $\Phi$ defined by $\Phi (g)=\obstruconj g$ and $\Phi (U_{i})=U_{i}$, for every $i\in\mathbb{Z}$. So $\Phi$ leaves $R_{-1}$ invariant and has order $p$. Let $\Phi _{k}=\Phi|_{R_{-1}^{\prime}\cap R_{k}}$. Then $\{\Phi _{k}\}_{k\geq 0}\in {\mathcal G}$. Recall that by Lemma \ref{lemma4.1}
\[R_{k}\cap R_{-1}^{\prime}=\{g, U_{-k}, U_{-k+1},\hdots,U_{-1}\}^{\prime\prime}.\]
Hence, by the previous remark, $\Phi _{k} (g)=\obstruconj g$ implies that $\Phi\not\in\chiusoRoR$.

Using the relations satisfied by the $U_{k}$ 's and $g$ one can easily verify that 
$$
\ad _{R_{0}}(U_{0})=\Phi\,\lim_{n\rightarrow\infty}\ad _{R_{0}}(U_{2}^{*}U_{4}^{*}\hdots U_{2n}^{*}),
$$
and   
$$
Ad_{R_{0}}(U_{-1})=\Phi\,\lim_{n\rightarrow\infty}Ad_{R_{0}}(U_{1}^{*}U_{3}^{*}\hdots U_{2n-1}^{*}).
$$
Note also that $\Phi$ commutes with the approximately inner automorphisms in the above decompositions. 

Therefore, the automorphism $\ad _{R_{0}}(U_{-1}^{k}U_{0}^{j})$ belongs to $\chiusoRoR$ if and only if $k+j\equiv 0\bmod p$.

To complete the proof we need to show that no non-trivial automorphism of the form $\displaystyle\ad _{R_{0}}\big (U_{-1}^{k}\sum _{i=1}^{p}\alpha _{i}g_{i}\big )$ is in $\chiusoRoR$. 
  
 Suppose $\displaystyle\ad _{R_{0}}\big (U_{-1}^{k}\sum _{i=1}^{p}\alpha _{i}g_{i}\big )\in\chiusoRoR$ for some choice of $0\leq k\leq p$ and $\alpha _{i}\in\mathbb{T}$, for all $1\leq i\leq p$. Since $\Ad _{R_{0}}(U_{-1})(g)=\obstruconj g$, we have that 

\begin{equation*}
\ad _{R_{0}}\left (U_{-1}^{k}\sum _{i=1}^{p}\alpha _{i}g_{i}\right )(g)=e^{-\frac{2\pi k i}{p}}g.
\end{equation*}

By the previous remark, it follows that $k=0$. Also, it is a consequence of the equality $\ad _{R_{0}}(U_{-1})(g_{k})=g_{k+1}$ for $1\leq k\leq p$, where $k+1$ is reduced $\bmod p$, that 
\[ \ad _{R_{0}}\left (\sum _{i=1}^{p}\alpha _{i}g_{i}\right )(U_{-1})=\left (\sum _{i=1}^{p}\alpha _{i}\bar{\alpha _{i-1}}g_{i}\right )U_{-1}.\]
Thus, $\ad _{R_{0}}(\sum _{k=1}^{p}\alpha _{k}g_{k})$ belongs to $\chiusoRoR$ if and only if $\alpha _{k}=\alpha _{k+1}$ for every $1\leq k< p$.

From the relations $\ad _{R_{0}}(U_{-1}^{p})=\ad _{R_{0}}(g)$, $U_{0}^{p^2}=1$ and $\ad _{R_{0}}(gU_{0}^{p})\in\chiusoRoR$, it follows that any automorphism of the form $\ad _{R_{0}}(U_{0}^{p-k}U_{-1}^{k})$ for $1\leq k\leq p^2$ can be written as a power of $\ad _{R_{0}}(U_{0}^{p-1}U_{-1})$, up to an element in $\innerRoR$. Indeed, if we choose $1\leq m\leq p-1$ such that $2m\equiv k-1\bmod p$ then 
\begin{eqnarray*}
&\hspace{-.5cm} &\ad _{R_{0}}(U_{0}^{p-k}U_{-1}^{k})=\ad _{R_{0}}(U_{0}^{p-k}U_{-1}^{-pm}U_{-1}^{k+pm})=\ad _{R_{0}}((g^{*})^{m}U_{0}^{p-k}U_{-1}^{k+pm})=\\
& \hspace{-.5cm}&  \ad _{R_{0}}((g^{*}U_{0}^{-p})^{m}U_{0}^{p+pm-k}U_{-1}^{k+pm})=\ad _{R_{0}}(g^{*}U_{0}^{-p})^{m}\ad _{R_{0}}(U_{0}^{p(1+2m)-pm-k}U_{-1}^{k+pm})=\\
&\hspace{-.5cm} & \ad _{R_{0}}(g^{*}U_{0}^{-p})^{m}\ad _{R_{0}}(U_{0}^{pk-pm-k}U_{-1}^{k+pm})=\ad _{R_{0}}(g^{*}U_{0}^{-p})^{m}\ad _{R_{0}}(U_{0}^{p-1}U_{-1})^{k+pm}.
\end{eqnarray*}
 \end{proof}

Using Corollary \ref{corollary5.5} we can now describe $\chi (\A )$    
\begin{corollary}
\label{corollary5.8}
Let $R_{-1}\subset R_{0}$ be an inclusion of type II$_{1}$ factors, where 
$R_{0}$ is the crossed product $R_{-1}\rtimes _{\ad_{R_{-1}}(U_{-1})}\luna$. 
Let $\A =\freeprodRRo$, with $Q$ a II$_{1}$ factor. Then $\chi (\A )\cong\luna$. Moreover, with 
the notation used before for the model $R_{-1}\subset R_{0}$, we have 
that 
$$
\chi (\A )=\langle\xi (Ad_{\A }(U_{0}^{p-1}U_{-1}))\rangle
$$
where $\xi$ denotes the quotient map from $\autA$ to $\displaystyle{\outA=\frac{\autA}{\innerA}}$.

\end{corollary}

Using the same type of argument used in \cite{Viola} and \cite{Jones3} one can use the exact sequence introduced by Connes in \cite{Connes1} to compute the Connes invariant of the crossed product $\crossprod$ since $\tenpro$ Has only trivial hypercentral sequences.

\begin{remark}
\label{remark5.9}
The II$_{1}$ factor $N=\tenpro$ has no non-trivial hypercentral sequence. Moreover any central sequence in $N$ is of the form $(1\otimes x_{n})_{n\geq 0}$, for a central sequence $(x_{n})_{n\geq 0}$ in $R_{0}$.
\end{remark}

 Indeed, since $\displaystyle N=\tenpro\cong ((*_{1} ^{n_{1}}Q)\otimes R_{0})*_{1\otimes R_{0}}((*_{1} ^{n_{2}}Q)\otimes R_{0})$ with $n_{1}+n_{2}=p$, we can apply Proposition \ref{proposition5.1} to $\tenpro$. If we denote by $E$ the conditional expectation on $1\otimes R_{0}$, and take any central sequence $(y_{n})_{n\geq 0}$ in $N$, then $E(y_{n})$ is a central sequence in $1\otimes R_{0}$ which approximate $(y_{n})$ in the $L^2$-norm.  

The proof of the following proposition follows the one given in \cite{Viola} for $p=3$, after observing that given any two II$_{1}$ factors $M_{1}$ and $M_{2}$ the free product factor $M_{1}*M_{2}$ is full (cf. \cite{Bar}).

\begin{proposition}
\label{proposition5.10}
Let $Q$ be a II$_{1}$ factor and $\M=\crossprod$. Then $\chi (\M)\cong\luna$.
\end{proposition}

\section{Non-outer conjugacy of the actions}

In the previous section we have shown that the $II_{1}$ factor $\A=\freeprodRRo$ has the same Connes invariant of the $II_{1}$ factor $\M=\crossprod$. However, the Connes invariant actually contains more information than what we have learn so far. This extra information derives from the position of this abelian group inside the group of outer automorphisms. To distinguish $\M$ and $\A$ we thus look at the crossed product of each of these two $II_{1}$ factors with the unique subgroup of order $p$ contained in $\chi (\M )\cong\chi (\A )\cong\luna$. We start this section by proving that the dual actions associated to this crossed products can be decomposed into an approximately inner automorphism and a centrally trivial automorphism. Using this decomposition we are able to show that $\A$ and $\M$ are not isomorphic. Hence, the $\luna$-actions defined in Theorem \ref{theorem4.8} and Corollary \ref{corollary4.10} are not outer conjugate. The argument we use here is a generalization to any choice of prime $p$ of the argument used by R\u{a}dulescu for $p=2$.

\begin{proposition}
\label{proposition6.1}
    Let $\A=\freeprodRRo$, with $Q$ a II$_{1}$ factor. The automorphism $s=Ad_{R_{0}}(U_{-1}^{p})$ defines an action of $\gruppo$ on $R_{0}$ which is the identity on $R_{-1}$. Moreover, the subalgebras ${\mathcal C}=\{U_{-1}^{p},\, gU_{0}^{p},\, U_{1}, U_{2}, \hdots\}^{\prime\prime}$ and ${\mathcal D}=\{U_{-1}^{p},\, g,\, U_{0},\, U_{1}, \hdots\}^{\prime\prime}$ of $R_{1}$ can be identified with the crossed products $\crossprodR$ and $\crossprodRo$, respectively. Both these von Neumann algebras are not factors: more precisely ${\mathcal Z}(\crossprodR )=\{U_{-1}^{p}\}^{\prime\prime}$ and  ${\mathcal Z}(\crossprodRo )=\{g^{*}U_{-1}^{p}\}^{\prime\prime}$.                

Lastly, 
$$
\A\rtimes _{Ad_{\A}(U_{-1}^{p})}\gruppo\cong (Q\otimes (\crossprodR))*_{\crossprodR}(\crossprodRo).
$$    

\end{proposition}
\begin{proof}
    Since $s (gU_{0}^{p})=gU_{0}^{p}$\, and $s (U_{i})=U_{i}$ for $i\geq 1$, we can conclude that $s$ acts as the 
    identity on $R_{-1}$. To show that  ${\mathcal D}$ can be identified with the crossed product $\crossprodRo$ we 
    need to verify that any monomial $m$ in ${\mathcal D}$ can be written using only one occurrence of $U_{-1}^{p}$ to     some power. Moreover, any monomial $m$ containing $U_{-1}^{p}$ must have zero trace in ${\mathcal D}$. 
   
    Both these properties are immediate consequences of the relations between the elements of $\{U_{k},g\mid k\in
    \mathbb{Z}\}$ and the definition of the trace on the von Neumann algebra generated by these 
    unitaries (see Lemma \ref{lemma4.1}). An analogous argument show that $\crossprodR=\{U_{-1}^{p},\, gU_{0}^{p},\, U_{1},\, U_{2},\hdots\}^{\prime\prime}$.
  
    In addition, from the observation that $s =Ad_{R_{0}}(U_{-1}^{p})$ acts identically on $R_{-1}$ it follows that 
    $U_{-1}^{p}$ belongs to the center of $\crossprodR$. Writing any element $x$ in $\crossprodR$ as $\displaystyle\sum _{j=0}^{p-1} x_{j}U_{-1} ^{pj}$ with $x_{j}\in R_{-1}$, it is easy to check that any element in the center of $\crossprodR$ belongs to $\{U_{-1}^{p}\}^{\prime\prime}$. Indeed, if $x$ commutes with every $y=\displaystyle\sum _{i=0}^{p-1} y_{i}U_{-1} ^{p i}\in \crossprodR$, then $x_{j}y_{i}=y_{i}x_{j}$ for all $1\leq i,j\leq p-1$, since $s$ acts identically on $R_{-1}$. Therefore, $y_{i}\in\mathbb{C}$ for all $1\leq i\leq p-1$. A similar argument is used to verify that $g^{*}U_{-1}^{p}$ 
    generates the center of $\crossprodRo$.  
    The claim that 
    $$
    \A\rtimes _{Ad_{\A}(U_{-1}^{p})}\gruppo\cong (Q\otimes (\crossprodR))*_{\crossprodR}(\crossprodRo)
    $$ 
    follows immediately from \cite[Remark 16]{Radu2}.  
\end{proof}

Let $\hat{s}$ be the dual action of $\widehat{\gruppo}$ on the crossed product $R_{0}\rtimes _{s}\gruppo$ defined above. Then $\hat{s}$ can be decomposed as a centrally trivial automorphism and an approximately inner automorphism. This proof is analogous to the one given by R\u{a}dulescu for the case $p=2$ and we present it here only for the convenience of the reader.

\begin{proposition}
\label{proposition6.2}
   Let $\hat {s}$ be the dual action of $\gruppo$ on $\crossprodRo=\{U_{-1}^{p},\, g,\, U_{0},\, U_{1}, \hdots\}^{\prime\prime}$ defined by $\hat{s} (x)=x$ for all $x\in R_{0}$, and $\hat{s} (U_{-1}^{p})=\obstruconj\, U_{-1}^{p}$. Set $\sigma=\ad _{\crossprodRo }(U_{-1}^{p-1}U_{0})$. Then, the decomposition $\hat{s}=(\hat{s}\sigma ^{-1})\sigma$ has the property that 
   \[ \hat{s}\sigma ^{-1}\in\chiuso,\]
   and
   \[\sigma\in\CtcrosRoR.\]
   Moreover, $\sigma\notin\inner$. 
\end{proposition}
\begin{proof}
   The relation $\ad _{\crossprodRo}(U_{-1})\mid _{\crossprodR }=Id\mid _{\crossprodR}$ implies immediately that
   $\ad _{\crossprodRo}(U_{-1})\in\CtcrosRoR$.
   
   Next we want to show that $\ad _{\crossprodRo}(U_{0})$ belongs to $\CtcrosRoR$. By Corollary \ref{corollary5.5}  we have that $\ad _{R_{0}}(U_{0})\in\CtRoR$. Since $R_{0}^{\omega}\cap R_{0}^{\prime}$ is the crossed product of $R_{-1}^{\omega}\cap R_{0}^{\prime}$ by $(\Ad U_{-1})_{\omega}$, it suffices to consider the central sequences for $\crossprodRo$ of the form $(U_{-1}^{pk}g^{m}U_{0}^{pm}U_{1}^{p(k+m)}U_{2}^{pm}U_{3}^{p(m+k)}\cdots U_{2n-1}^{p(m+k)}U_{2n}^{pk})_{n\in\mathbb{N}}$, for $1\leq k,n\leq p$. Thus, it is enough to see how $\Ad U_{0}$ acts on $(U_{-1}^{p}U_{1}^{p}U_{3}^{p}\cdots U_{2n-1}^{p})_{n\in\mathbb{N}}$, since $(gU_{0}^{p}U_{1}^{p}U_{2}^{p}U_{3}^{p}\cdots U_{2n-1}^{p}U_{2n}^{p})_{n\in\mathbb{N}}$ is a central sequence of $\crossprodRo$ which belongs to $R_{-1}$ and $\ad U_{0}$ acts trivially on such sequences. A simple computation shows that  $\Ad U_{0}$ leaves $(U_{-1}^{p}U_{1}^{p}U_{3}^{p}\cdots U_{2n-1}^{p})_{n\in\mathbb{N}}$ invariant. Therefore, $\Ad U_{0}$ belongs to $\CtcrosRoR$, as well as $\sigma =Ad_{\crossprodRo}(U_{-1}^{p-1}U_{0})$.
   
   Also, using the properties satisfied by the $U_{k}$'s and $g$ it is easy to show that if
   \begin{equation*}
   x_n=\begin{cases}
       g^{*} U_{0}^{-p}U_{1}^{*}U_{2}^{1-p}U_{3}^{*}U_{4}^{1-p}\hdots U_{n}^{*}, 
       &\text{ if $n$ is odd}, \\
       g^{*} U_{0}^{-p}U_{1}^{*}U_{2}^{1-p}U_{3}^{*}U_{4}^{1-p}\hdots U_{n}^{1-p},& 
       \text{ if $n$ is even}. \\
       \end{cases} 
   \end{equation*}  
   then
\[\sigma ^{-1}=\lim_{n\rightarrow\infty}\ad _{R_{0}}(x_{n})\in\chiusoRoR.\]

   However, $\sigma ^{-1}\not\in\chiuso$.
   Indeed,
   \begin{equation*}
   \sigma ^{-1}(g^{*}U_{-1}^{p})=\obstruction\, g^{*}U_{-1}^{p},
   \end{equation*}
   while $g^{*}U_{-1}^{p}$ is left invariant by the inner automorphisms of 
   $\crossprodRo$ since it belongs to the center of $\crossprodRo$ (Proposition \ref{proposition6.1}). However, multiplying $\sigma ^{-1}$ by $\hat{s}$ we obtain that 
   \begin{equation*}
   \hat{s}\sigma ^{-1}=\lim_{n\rightarrow\infty}\ad _{\crossprodRo}(x_{n})\in \chiuso.
   \end{equation*}

   In addition, since $\sigma (U_{-1}^{p})=\obstruction U_{-1}^{p}$ and $U_{-1}^{p}\in{\mathcal Z}(\crossprodR )$, it follows that $\sigma$ does not belong to $\operatorname{Int}(\crossprodRo , \crossprodR )$.
\end{proof}

Observe that the map $\Phi :\operatorname{Aut}(\crossprodRo, \crossprodR)\longrightarrow\autB$ defined by $\displaystyle\Phi (\alpha)=(Id\otimes\alpha |_{(\crossprodR}))*_{(\crossprodR)}\alpha$ is continuous and sends $\operatorname{Int}(\crossprodRo,\crossprodR )$ into $\innerB$, and $\operatorname{Ct}(\crossprodRo,\crossprodR)$ into $\CtB$ \cite[Proposition 3]{Radu2}. Thus, the previous decomposition of the dual action on $R\rtimes _{s}\gruppo$ yields a decomposition of the dual action on the crossed product $\B=\A\rtimes _{\ad _{\A}(U_{-1}^{p})}\gruppo$ into an approximately inner automorphism and a centrally trivial automorphism.

\begin{corollary}
\label{corollary6.3}
   Let $R_{-1}\subset R_{0}$ be the inclusion of type $II _{1}$ factors
   of lemma \ref{lemma4.1}. Set $\A =\freeprodRRo$, with $Q$ a II$_{1}$ factor. Let 
   $\B =\A\rtimes _{\ad _{\A}(U_{-1}^{p})}\gruppo$, which 
   by Proposition \ref{proposition6.1} is isomorphic to  
   $(Q\otimes (\crossprodR))*_{\crossprodR}(\crossprodRo)$. Since the $\gruppo$ action realizing the crossed product generates the only subgroup of order $p$ in $\chi (\A)\cong\luna$, it is an invariant of the factor $\A$. Denote 
   by $S$ the dual action on the crossed product $\B$. 

   \noindent Then $S$ can be decomposed as
   $$
   S=S_{1}S_{2},
   $$
   where $S_{1}=S\ad _{\B}(U_{-1}^{p-1}U_{0})^{*}\in\chiusoB$ and $S_{2}=
   \ad _{\B}(U_{-1}^{p-1}U_{0})\in\CtB$. Moreover, if we set $h=gU_{-1}^{p}U_{0}^{p}$, then 
   \begin{equation}
   \label{eqn6.3}
   \begin{array}{l}
   S_{1}^{p}=\ad _{\B}h^{*}\\
   S_{2}^{p}=\ad _{\B}h,
   \end{array}\end{equation}
   and 
   $$
   S_{i}(h)=\obstruction h, \text{ for }i=1,\, 2.
   $$

\end{corollary}

Now consider the crossed product $\M=\crossprod$. Using Takesaki duality we can easily get that the dual action on the crossed product $\dualcr$ has a unique decomposition into an approximately inner automorphism and a centrally trivial automorphism, similar to the one found for the dual action $S$ on the factor $\B=\A\rtimes _{\ad _{\A}(U_{-1}^{p})}\gruppo$. The uniqueness of the decomposition is a consequence of the following remark.

\begin{remark}
\label{remark6.4} 
The II$_{1}$ factor $P=\tenpro$ has trivial Connes invariant. 
\end{remark}   
Indeed, assume $\alpha$ is an automorphism in $\CtP\cap\chiusoP$. The assumption $\alpha\in\CtP$ implies that $\alpha=\Ad z(\nu\otimes Id)$ for some unitary $z\in\tenpro$ and an automorphism $\nu$ of $\ifreeprod$ (simply replace $\mathfrak{L}(\Fb  _{t})$ with $\ifreeprod$ in \cite[Lemma3.6]{Viola}). Since $\alpha\in\chiusoP$, we obtain that $\nu$ belongs to $\overline{\operatorname{Int}\left (\ifreeprod\right )}=\operatorname{Int}\left (\ifreeprod\right )$.
 
\begin{proposition}
\label{proposition6.5}
   Let $P=\tenpro$ and $\M=P\rtimes _{\gamma}\gruppo$. The dual action $\widetilde{\gamma}$ on the crossed product $\dualcr$ can be uniquely decomposed, up to an inner automorphism, as 
   $\widetilde{\gamma}=\Ad w\widetilde{\gamma} _{1}\widetilde{\gamma} _{2}$, where
   $\widetilde{\gamma} _{1}\in\overline{\operatorname{Int}(\dualcr)}$, 
   $\widetilde{\gamma} _{2}\in\operatorname{Ct}(\dualcr)$, and $w$ is a 
   unitary in $\dualcr$ . 

   In addition, $\widetilde{\gamma} _{1}$ and $\widetilde{\gamma} _{2}$ have outer period $p$ and conjugate obstruction to lifting. More precisely, there exist unitaries $f_{i}$ for $i=1, 2$, satisfying the 
   conditions
   \[\begin{array}{l}
             \widetilde{\gamma} _{i}^{p}=\ad _{\dualcr}(f_{i}), 
             \widetilde{\gamma} _{1}(f_{1})=\obstruconj f_{1}, \\
             \widetilde{\gamma} _{2}(f_{2})=\obstruction f_{2},
     \end{array}\]
   and
   \[\widetilde{\gamma} _{i}(f_{j})=f_{j},\text{ for }i\neq j.\]

\end{proposition}
\begin{proof}
     By Takesaki duality \cite[Theorem 4.6]{Take} 
     $$
     \dualcr\cong P\otimes B(\ell ^{2}(\gruppo)).
     $$
     Moreover, the dual action $\widetilde{\gamma}$ of $\widehat{\gamma}$ corresponds under this identification with the action $\gamma\otimes \operatorname{Ad}\lambda(1)^{*}$, where $\lambda$ is the usual left representation of $\gruppo$ on $\ell ^{2}(\gruppo )$ defined by 
\[(\lambda (h)\eta)(k)=\eta (k-h), \text{ for }h,k\in\gruppo,\;\eta\in\ell ^{2}(\gruppo).\]
     Observe that 
     $$
     \gamma =\Ad W\,(\alpha\otimes\beta )=\Ad W\,(1\otimes\beta )
     (\alpha\otimes 1),
     $$
     and recall that by Remark \ref{remark5.9} any central sequence in $P$ has the form 
     $(1\otimes y_{n})$, for a central sequence $(y_{n})$ in $R_{0}$. Therefore, $\alpha\otimes 1\in\CtN$. Furthermore, the sequence ($x_{n})_{n\geq 0}$ of unitaries in $P$ given by
      \begin{equation*}
       x_n=\begin{cases}
        U_{0}U_{1}^{*}U_{2}U_{3}^{*}\hdots U_{n}^{*},&\text{ if $n$ is odd}, \\
        U_{0}U_{1}^{*}U_{2}U_{3}^{*}\hdots U_{n},& \text{ if $n$ is even} \\
       \end{cases} 
     \end{equation*} 
     has the property $\displaystyle{\beta=\lim_{n\rightarrow\infty}\Ad_{R_{0}}     (x_{n})}$. Thus $1\otimes\beta\in\chiusoP$. 
         
     Setting $\gamma _{1}=1\otimes\beta$ and $\gamma _{2}=\alpha\otimes 1$,
     we can easily check that the following relations are satisfied:
     \begin{equation*}
     \gamma _{1}^{p}=Ad\, (1\otimes g)\text{ with }
     \gamma _{1}(1\otimes g)=\obstruconj (1\otimes g),
     \end{equation*}
     and
     \begin{equation*}
     \gamma _{2}^{p}=Ad\, (u\otimes 1)\text{ with }
     \gamma _{2}(u\otimes 1)=\obstruction (u\otimes 1).
     \end{equation*}
     Denote by $Id$ the identity of $B(\ell ^{2}(\gruppo))$. Let $\widetilde{\gamma} _{1}=\gamma _{1}\otimes\operatorname{Ad}(\lambda(1 ))^{*}$ and $\widetilde{\gamma} _{2}=\gamma _{2}\otimes Id$. Then, $\widetilde{\gamma} _{1}\in\overline{\operatorname{Int}(\dualcr)}$ and $\widetilde{\gamma} _{2}\in\operatorname
    {Ct}(\dualcr)$. Set $f_{1}=(1\otimes g)\otimes Id$, $f_{2}=(u\otimes 1)\otimes Id_{B(\ell ^{2}(\gruppo))}$ and $w=W\otimes Id$. Obviously,  
    \begin{equation*}
        \widetilde{\gamma} _{1}^{p}=\ad f_{1}\,\text{ with }\,
         \widetilde{\gamma} _{1}(f_{1})=\obstruconj f_{1},
     \end{equation*}
     and
     \begin{equation*}
        \widetilde{\gamma} _{2}^{p}=\ad f_{2}\,\text{ with }\,
        \widetilde{\gamma} _{2}(f_{2})=\obstruction f_{2}.
     \end{equation*}
 
     \noindent Moreover, $\widetilde{\gamma} _{1}(f_{2})=f_{2}$, $\widetilde{\gamma}
      _{2}(f_{1})=f_{1}$ and $\widetilde{\gamma}=\Ad w\widetilde{\gamma}
      _{1}\widetilde{\gamma} _{2}$. 
\end{proof}

Using the decompositions of the dual actions $S$ and $\widehat{\gamma}$, we can now prove our main theorem which implies the non outer-conjugacy of the two $\luna$-actions defined on $\freeprodN$. Several arguments used in the proof are due to R\u{a}dulescu.  

\begin{theorem}
\label{theorem6.6}
      Let $R_{-1}\subset R_{0}$ be an inclusion of type $II _{1}$ factors, 
    where $R_{0}$ is the crossed product $R_{-1}\rtimes _{\theta}\luna$ and $\theta$ has outer invariant $(p,\obstruconj)$. Given a II$_{1}$ factor Q, let $\A=\freeprodRRo$ and $\M=\crossprod$. Then $\M$ is not 
    isomorphic to $\A$.
\end{theorem}
\begin{proof}
    We will prove it by contradiction. Assume that $\M$ and $\A$ are isomorphic. Recall that $\widehat{\gamma}$ and $\Ad _{\A}(U_{-1}^{p})$, respectively, generate the only subgroup of order $p$ in $\chi (\M)\cong\luna$ and $\chi (\A)\cong\luna$, respectively. Therefore, $\M _{1}=\dualcr$ and $\B =\A\rtimes _{\ad _{\A}(U_{-1}^{p})}\gruppo$ are also isomorphic, and the dual actions defined on these crossed products differ only by an inner automorphism. It follows that the decomposition of $\widetilde{\gamma}$ in a centrally trivial automorphism and an approximately inner automorphism (Proposition \ref{proposition6.5}) must be outer conjugate to the similar decomposition given for $S$ in Corollary \ref{corollary6.3}. 

Denote 
    by $\Gamma$ the isomorphism between $\B=\A\rtimes _{\ad _{\A}(U_{-1}^{p})}\gruppo$ and 
    $\M _{1}=\M\rtimes _{\hat{\gamma}}\gruppo$, and set $S_{1}^{\prime}=\Gamma 
    S_{1}\Gamma ^{-1}$, $S_{2}^{\prime}=\Gamma S_{2}\Gamma ^{-1}$, where $S_{1}$ and $S_{2}$ are the two automorphisms appearing in the decomposition of $S$ in Corollary \ref{corollary6.3}. Then, for some 
    unitary $v_{0}$ we have that
    \begin{equation*}
    S_{1}^{\prime}S_{2}^{\prime}=\Ad v_{0}\widetilde{\gamma} _{1}\widetilde{\gamma} _{2}.
    \end{equation*}
    Set $S^{\prime}=S_{1}^{\prime}S_{2}^{\prime}$ and $T=\widetilde{\gamma} _{1}\widetilde{\gamma}
     _{2}$, so that $S^{\prime}=\Ad (v_{0})T$.
    Observe that $S^{\prime}$ has period $p$ while $T$ has period $p^2$, so that  
    \begin{equation*}
    1=(S^{\prime})^{p^2}=\left (\Ad \prod _{j=0} ^{p^{2}-1}{T^{j}(v_{0})}\right ).
    \end{equation*}
    Since $\dualcr$ is a factor, this implies that 
    \begin{equation*}
    \prod _{j=0} ^{p^{2}-1}{T^{j}(v_{0})}=\mu\, 1,
    \end{equation*}
    for some $\mu\in\mathbb{T}$. 

    Let $\mu '$ be a $p^{2}$-root of $\mu$. Set $v_{0}'=\frac{1}{\mu '}v_{0}$. Then 
    $v_{0}'$ satisfies the relation
    \begin{equation*}
     \prod _{j=0} ^{p^{2}-1}{T^{j}(v_{0}')}=1.
    \end{equation*}
    By part (i) of Corollary 2.6 in \cite{Connes2} we obtain that $v_{0}'$
    has the form $w^{*}T(w)$ for some unitary $w$ in $\dualcr$. Set 
    $$
    \Gamma _{1}=\Ad (w^{*})\widetilde{\gamma} _{1}\Ad (w)\quad\text{ and }
    \quad\Gamma _{2}=\Ad (w^{*})\widetilde{\gamma} _{2}\Ad (w).    
    $$ 
    Then $\Gamma _{1}$ and $\Gamma _{2}$ commute, and
    \begin{equation}
    \label{arancio}
    S_{1}^{\prime}S_{2}^{\prime}=\Gamma _{1}\,\Gamma _{2},
    \end{equation}
    where $S_{1}^{\prime}, \Gamma _{1}\in\overline{\operatorname{Int}(\dualcr)}$ 
    and
    $S_{2}^{\prime}, \Gamma _{2}\in\operatorname{Ct}(\dualcr)$. Therefore $S_{2}^{\prime}\Gamma _{2}^{-1}=(S_{1}^{\prime})^{-1}\Gamma _{1}\in\operatorname{Ct}(\dualcr)\cap\overline{\operatorname{Int}(\dualcr)}$.  
    
    Takesaki duality and  Remark \ref{remark6.4} imply that   
    \begin{equation*}
    S_{2}^{\prime}=Ad\, w_{2}\,\Gamma _{2}\quad\text{ and }\quad
    S_{1}^{\prime}=\Gamma _{1}Ad\, w_{2}^{*}=Ad\, w_{1}\,\Gamma _{1},
    \end{equation*}
    for some unitary $w_{2}\in \dualcr$ and $w_{1}=\Gamma _{1}(w_{2}^{*})$.
    Moreover, for $i,j\in\{1,2\}$, $i\neq j$ we have
    \begin{equation} 
    \label{mela1}
    \Gamma _{i}(w^{*}f_{j}w)=w^{*}f_{j}w.
    \end{equation}
    \begin{equation}
    \label{mela}
    (\Gamma _{i})^{p}=\Ad (w^{*}f_{i}w), 
    \end{equation}
    
    \noindent In addition,
    \begin{equation}
    \label{pesca}
    \Gamma _{1}(w^{*}f_{1}w)=\obstruconj\, w^{*}f_{1}w\quad\text {and}\quad 
    \Gamma _{2}(w^{*}f_{2}w)=\obstruction\, w^{*}f_{2}w.
    \end{equation}
    
    Using (\ref{mela}) and (\ref{pesca}) we also obtain
    \begin{align}
    & (S_{i}^{\prime})^{p}=(\Ad w_{i}\Gamma _{i})^{p}=\Ad \left (\left (\prod _{j=0} ^{p-1}{\Gamma _{i} ^{j}(w_{i})}\right )w^{*}f_{i}w\right ), 
    \text{ for }i=1,2\\
    & S_{1}^{\prime}\left (\left (\prod _{j=0} ^{p-1}{\Gamma _{1} ^{j}(w_{1})}\right )w^{*}f_{1}w\right )=\obstruconj\, 
    \left (\prod _{j=0} ^{p-1}{\Gamma _{1} ^{j}(w_{1})}\right )w^{*}f_{1}w, \\
    & S_{2}^{\prime}\left (\left (\prod _{j=0} ^{p-1}{\Gamma _{2} ^{j}(w_{2})}\right )w^{*}f_{2}w\right )=\obstruction\, \left (\prod _{j=0} ^{p-1}{\Gamma _{2} ^{j}(w_{2})}\right )w^{*}f_{2}w.
    \end{align}
   
    Since $[S_{1}^{\prime},S_{2}^{\prime}]=0$ and $[\Gamma _{1},\Gamma _{2}]=0$, we have that
    $$
    \Ad (w_{1}\Gamma _{1}(w_{2}))\Gamma _{1}\Gamma _{2})=S_{1}^{\prime}S_{2}^{\prime}=S_{2}^{\prime}S_{1}^{\prime}=
    \Ad (w_{2}\Gamma _{2}(w_{1}))\Gamma _{1}\Gamma _{2},
    $$
    and $\Gamma _{1}$ and $\Gamma _{2}$ commute, we conclude 
    that
    $$
    \Ad (w_{1}\Gamma _{1}(w_{2}))=\Ad (w_{2}\Gamma _{2}(w_{1})), 
    $$
    which implies that there exists a complex number $\lambda$ of modulus 1
    such that
    \begin{equation}
    \label{ban}
    w_{1}\Gamma _{1}(w_{2})=\lambda\, w_{2}\Gamma _{2}(w_{1}).
    \end{equation}
    
    Next we use (\ref{mela1}), (\ref{mela}), (\ref{ban}) and the fact that $\Gamma _{1}$ and $\Gamma _{2}$ commute to evaluate
    $\displaystyle S_{1}^{\prime}\left (\left (\prod _{j=0} ^{p-1}{\Gamma _{2} ^{j}(w_{2})}\right )w^{*}f_{2}w\right )$.
    
    \begin{equation*}
     \begin{split}
        & S_{1}^{\prime}\left (\left (\prod _{j=0} ^{p-1}{\Gamma _{2} ^{j}(w_{2})}\right )w^{*}f_{2}w\right )=
        \Ad w_{1}\Gamma _{1}\left (\left (\prod _{j=0} ^{p-1}{\Gamma _{2} ^{j}(w_{2})}\right )w^{*}f_{2}w\right )= \\
        & w_{1}\left (\prod _{j=0} ^{p-1}{\Gamma _{1}\Gamma _{2} ^{j}(w_{2})}\right )\Gamma _{1}(w^{*}f_{2}w )w_{1}^{*}=\lambda 
        w_{2}\Gamma _{2}(w_{1}\Gamma _{1}(w_{2}))\Gamma _{1}\Gamma _{2}^{2}(w_{2})\cdots\Gamma _{1}\Gamma _{2}^{p-1}(w_{2})w^{*}f_{2}ww_{1}^{*}= \\
        & \lambda ^{2}w_{2}\Gamma _{2}(w_{2})\Gamma _{2}^{2}(w_{1}\Gamma _{1}
        (w_{2}))\cdots\Gamma _{1}\Gamma _{2}^{p-1}(w_{2})w^{*}f_{2}ww_{1}^{*}=\lambda ^{p}\left (\prod _{j=0} ^{p-1}{\Gamma _{2} ^{j}(w_{2})}\right )\Gamma _{2} ^{p}(w_{1})w^{*}f_{2}ww_{1}^{*}=\\
& \lambda^{p} \left (\prod _{j=0} ^{p-1}{\Gamma _{2} ^{j}(w_{2})}\right )w^{*}f_{2}w.
    \end{split}
    \end{equation*}

    Because of the relations (\ref{eqn6.3}), and the fact that $\dualcr$ is a factor, we have that 
$$
\left (\prod _{j=0} ^{p-1}{\Gamma _{2} ^{j}(w_{2})}\right )w^{*}f_{2}w=\delta \left (\prod _{j=0} ^{p-1}{\Gamma _{1} ^{j}(w_{1})}\right )w^{*}f_{1}w)^{*}, 
$$
for some complex $\delta$ of modulus one.    Using (6) we conclude that 
\begin{equation}
\label{eqn6.5}
\lambda ^{p}=\obstruction.
\end{equation} 

    Lastly, using the relations $S_{i}^{\prime}=\Ad w_{i}\,\Gamma _{i}$ for 
    $i=1,2$, we obtain the equation 
    $$
    \Gamma _{1}\Gamma _{2}= S_{1}'S_{2}'=(\Ad w_{1}\,\Gamma _{1})(\Ad w_{2}\,\Gamma _{2})=\Ad (w_{1}\Gamma _{1}(w_{2}))\Gamma _{1}\Gamma _{2},
    $$
    which implies that 
    \[\Ad (w_{1}\Gamma _{1}(w_{2}))=Id.\]
    Hence, there exists a complex number $\lambda _{1}\in\mathbb{T}$ such that
    \begin{equation*}
    \Gamma _{1}(w_{2})=\lambda _{1}w_{1}^{*}.
    \end{equation*}
    Analogously, from $S_{2}'S_{1}'=\Gamma _{2}\Gamma _{1}$
    it follows that there exists $\lambda _{2}\in\mathbb{T}$ such that
    \begin{equation*}
    \Gamma _{2}(w_{1})=\lambda _{2}w_{2}^{*}.
    \end{equation*}
    These two relations, together with (\ref{ban}), imply that 
    $\lambda=\frac{\lambda _{1}}
    {\lambda _{2}}$. To get a contradiction to the original assumption that $\M\cong\A$ we evaluate $\Gamma _{1}
    \Gamma _{2}(w_{1})$:
    \begin{equation*}
    \Gamma _{1}\Gamma _{2}(w_{1})=\lambda _{2}\Gamma _{1}(w_{2}^{*})=
    \frac{\lambda _{2}}{\lambda _{1}}w_{1}=\frac{1}{\lambda}w_{1}.
    \end{equation*}
    On the other hand, $S_{1}'S_{2}'=\Gamma _{1}\Gamma _{2}$ has period $p$, so  $\lambda$ must be a $p$-root of unity, contradicting (\ref{eqn6.5}).
\end{proof}
    
\begin{corollary}
\label{corollary6.7}
Given any prime $p$ and any II$_{1}$ factor $N$ there exists two $\luna$-actions on $\freeprodN$ which have the same outer invariant $(p,\obstruction)$ but are not outer conjugate.
\end{corollary}
\begin{proof}
The $\luna$-actions $\widetilde{\theta _{1}}$ and $\widetilde{\theta _{2}}$ defined in Theorem \ref{theorem4.8} and Corollary \ref{corollary4.10} have the desired outer invariant. If they were outer conjugate, then the actions $\theta_{1}$ (Theorem \ref{theorem4.8}) and $\theta_{2}$ (Theorem \ref{theorem4.9}) would also be outer conjugate. Take the II$_{1}$ factor $Q$ so that $Q_{t}=N\otimes M(\mathbb{C})$, for $t=\sqrt{\frac{p+1}{p^3}}$. Observe that if $\theta_{1}$ and $\theta_{2}$ were outer conjugate, then $\M=\crossprod$ and $\A=\freeprodRRo$ would be isomorphic. This is because $\M$ and $\A$ are the enveloping algebras of a subfactor construction of the form $A_{i}\subset B_{i}=A_{i}\rtimes _{\theta _{i}}\luna$ for $i=1,2$, as described in Proposition \ref{proposition4.7} and Theorem \ref{theorem4.9}. But by our previous result such an isomorphism cannot exist.
\end{proof}
 \section*{Acknowledgments}

The second named author thanks Professor T.\ Giordano, Professor Y.\ Kawahigashi and Professor V.\ F.\ Jones for many useful discussions regarding this manuscript.

\bibliographystyle{plain}
            
\end{document}